\documentclass[a4paper,10pt]{article}
\usepackage[english]{babel}
\usepackage[latin1]{inputenc}
\usepackage{amssymb}
\usepackage{amsmath}
\usepackage{amscd}
\usepackage{latexsym}
\usepackage{amsthm}
\usepackage{eucal}
\usepackage{mathrsfs}
\theoremstyle{plain}
\newtheorem{thm}{Theorem}[section]
\newtheorem{cor}[thm]{Corollary}
\newtheorem{lem}[thm]{Lemma}
\newtheorem{prop}[thm]{Proposition}

\newtheorem{defn}[thm]{Definition}
\newtheorem{oss}[thm]{Remark}
\newtheorem{esem}[thm]{Example}

\title{Deformation of the O'Grady moduli spaces}
\author{Arvid Perego, Antonio Rapagnetta}

\begin{document}

\maketitle

\begin{abstract}
In this paper we study moduli spaces of sheaves on an abelian or projective K3 surface. If $S$ is a K3, $v=2w$ is a Mukai vector on $S$, where $w$ is primitive and $w^{2}=2$, and $H$ is a $v-$generic polarization on $S$, then the moduli space $M_{v}$ of $H-$semistable sheaves on $S$ whose Mukai vector is $v$ admits a symplectic resolution $\widetilde{M}_{v}$. A particular case is the $10-$dimensional O'Grady example $\widetilde{M}_{10}$ of irreducible symplectic manifold. We show that $\widetilde{M}_{v}$ is an irreducible symplectic manifold which is deformation equivalent to $\widetilde{M}_{10}$ and that $H^{2}(M_{v},\mathbb{Z})$ is Hodge isometric to the sublattice $v^{\perp}$ of the Mukai lattice of $S$. Similar results are shown when $S$ is an abelian surface.
\end{abstract}

\section{Introduction and notations}

Moduli spaces of semistable sheaves on abelian or projective K3 surfaces are an important tool to produce examples of irreducible symplectic manifolds. In the following, $S$ will denote an abelian or projective K3 surface. 

An element $v\in \widetilde{H}(S,\mathbb{Z}):=H^{2*}(S,\mathbb{Z})$ will be written it as $v=(v_{0},v_{1},v_{2})$, where $v_{i}\in H^{2i}(S,\mathbb{Z})$, and $v_{0},v_{2}\in\mathbb{Z}$. If $v_{0}\geq 0$ and $v_{1}\in NS(S)$, then $v$ is called \textit{Mukai vector}. Recall that $\widetilde{H}(S,\mathbb{Z})$ has a pure weight-two Hodge structure defined as $$\widetilde{H}^{2,0}(S):=H^{2,0}(S),\,\,\,\,\,\,\,\,\,\,\,\widetilde{H}^{0,2}(S):=H^{0,2}(S),$$ $$\widetilde{H}^{1,1}(S):=H^{0}(S,\mathbb{C})\oplus H^{1,1}(S)\oplus H^{4}(S,\mathbb{C}),$$and a lattice structure with respect to the Mukai pairing $(.,.)$. In the following, we let $v^{2}:=(v,v)$ for every Mukai vector $v$; moreover, for every Mukai vector $v$ define the sublattice $$v^{\perp}:=\{\alpha\in \widetilde{H}(S,\mathbb{Z})\,|\,(\alpha,v)=0\}\subseteq\widetilde{H}(S,\mathbb{Z}),$$which inherits a pure weight-two Hodge structure from the one on $\widetilde{H}(S,\mathbb{Z})$.

If $\mathscr{F}$ is a coherent sheaf on $S$, we define its \textit{Mukai vector} to be $$v(\mathscr{F}):=ch(\mathscr{F})\sqrt{td(S)}=(rk(\mathscr{F}),c_{1}(\mathscr{F}),ch_{2}(\mathscr{F})+\epsilon rk(\mathscr{F})),$$where $\epsilon=1$ if $S$ is K3, and $\epsilon=0$ if $S$ is abelian. Let $H$ be an ample line bundle on $S$. For every $n\in\mathbb{Z}$ and every coherent sheaf $\mathscr{F}$, let $\mathscr{F}(nH):=\mathscr{F}\otimes\mathscr{O}_{S}(nH)$. The Hilbert polynomial of $\mathscr{F}$ with respect to $H$ is $P_{H}(\mathscr{F})(n):=\chi(\mathscr{F}(nH))$, and the reduced Hilbert polynomial of $\mathscr{F}$ with respect to $H$ is $$p_{H}(\mathscr{F}):=\frac{P_{H}(\mathscr{F})}{\alpha_{H}(\mathscr{F})},$$where $\alpha_{H}(\mathscr{F})$ is the coefficient of the term of highest degree in $P_{H}(\mathscr{F})$.

\begin{defn}
A coherent sheaf $\mathscr{F}$ is $H-$\textit{stable} (resp. $H-$\textit{semistable}) if it is pure and for every proper $\mathscr{E}\subseteq\mathscr{F}$ we have $p_{H}(\mathscr{E})(n)<p_{H}(\mathscr{F})(n)$ (resp. $p_{H}(\mathscr{E})(n)\leq p_{H}(\mathscr{F})(n)$) for $n\gg 0$. 
\end{defn}

Let $H$ be a polarization and $v$ a Mukai vector on $S$. We write $M_{v}(S,H)$ (resp. $M_{v}^{s}(S,H)$) for the moduli space of $H-$semistable (resp. $H-$stable) sheaves on $S$ with Mukai vector $v$. If no confusion on $S$ and $H$ is possible, we drop them from the notation. 

From now on, we suppose that $H$ is a $v-$generic polarization (for a definition, see section 2.1). We write $v=mw$, where $m\in\mathbb{N}$ and $w$ is a primitive Mukai vector on $S$. It is known that if $M_{v}^{s}\neq\emptyset$, then $M_{v}^{s}$ is smooth, quasi-projective, of dimension $v^{2}+2$ and carries a symplectic form (see Mukai \cite{M1}). If $S$ is abelian, a further construction is necessary: choose $\mathscr{F}_{0}\in M_{v}(S,H)$, and define $a_{v}:M_{v}(S,H)\longrightarrow S\times\widehat{S}$ in the following way (see \cite{Y4}): let $p_{\widehat{S}}:S\times \widehat{S}\longrightarrow\widehat{S}$ be the projection and $\mathscr{P}$ the Poincar\'e bundle on $S\times\widehat{S}$. For every $\mathscr{F}\in M_{v}(S,H)$ we let $$a_{v}(\mathscr{F}):=(det(p_{\widehat{S}!}((\mathscr{F}-\mathscr{F}_{0})\otimes(\mathscr{P}-\mathscr{O}_{S\times\widehat{S}})),det(\mathscr{F})\otimes det(\mathscr{F}_{0})^{-1}).$$Moreover, we define $K_{v}(S,H):=a_{v}^{-1}(0_{S},\mathscr{O}_{S})$, where $0_{S}$ is the zero of $S$. Moreover, we denote $\widehat{H}$ the dual polarization on $\widehat{S}$.

If $S$ is a projective K3 surface, $v=w$ (i. e. $m=1$) and $H$ is $v-$generic, the moduli space $M_{v}(S,H)$ is well understood thanks to the work of several authors (see Mukai \cite{M2}, Beauville \cite{B}, O'Grady \cite{OG1}, Yoshioka \cite{Y1}, \cite{Y2}). The final result is the following:

\begin{thm}
\label{thm:vprimitivo}
(\textbf{Mukai, Yoshioka}). Let $S$ be an abelian or projective K3 surface, $v$ a primitive Mukai vector and $H$ a $v-$generic polarization. Then $M_{v}(S,H)=M_{v}^{s}(S,H)$, and we have the following results:
\begin{enumerate}
 \item if $S$ is K3 and $v^{2}=-2$, then $M_{v}(S,H)$ is a reduced point; 
 \item if $S$ is K3 and $v^{2}=0$, then $M_{v}(S,H)$ is a K3 surface and there is a Hodge isometry between $v^{\perp}/\mathbb{Z}\cdot v$ and $H^{2}(M_{v},\mathbb{Z})$; 
 \item if $S$ is K3 and $v^{2}\geq 2$, then $M_{v}$ is an irreducible symplectic variety of dimension $2n=v^{2}+2$, which is deformation equivalent to $Hilb^{n}(S)$, the Hilbert scheme of $n-$points on $S$. Moreover, there is a Hodge isometry between $v^{\perp}$ and $H^{2}(M_{v},\mathbb{Z})$, where the latter has a lattice structure given by the Beauville form;
 \item if $S$ is abelian and $v^{2}\geq 6$, then $K_{v}(S,H)$ is an irreducible symplectic variety of dimension $2n=v^{2}-2$, which is deformation equivalent to $K^{n}(S)$, the generalized Kummer variety on $S$, and there is a Hodge isometry between $v^{\perp}$ and $H^{2}(K_{v},\mathbb{Z})$.
\end{enumerate}
\end{thm}

If $v$ is not primitive, then $M_{v}$ can be singular: it is natural to ask if there is a \textit{symplectic resolution} of $M_{v}$, i. e. a resolution of singularities $\pi_{v}:\widetilde{M}_{v}\longrightarrow M_{v}$ such that on $\widetilde{M}_{v}$ there is a symplectic form extending the one on $M_{v}^{s}$. The first result appearing in the literature is the following:

\begin{thm}
\label{thm:og}(\textbf{O'Grady}).
\begin{enumerate}
 \item Let $S$ be a projective K3 surface, $v:=(2,0,-2)$ and $H$ a $v-$generic polarization. Then $M_{10}:=M_{v}(S,H)$ admits a symplectic resolution $\pi:\widetilde{M}_{10}\longrightarrow M_{10}$, and $\widetilde{M}_{10}$ is an irreducible symplectic variety of dimension 10 and second Betti number 24.
 \item Let $S$ be an abelian surface, $v:=(2,0,-2)$ and $H$ a $v-$generic polarization. Then $K_{6}:=K_{v}(S,H)$ admits a symplectic resolution $\pi:\widetilde{K}_{6}\longrightarrow K_{6}$, and $\widetilde{K}_{6}$ is an irreducible symplectic variety of dimension 6 and second Betti number 8.
\end{enumerate}
\end{thm}

The first example is studied by O'Grady in \cite{OG2}, the computation of its second Betti number is started in \cite{OG2} and completed in \cite{R1}. The second example is contained in \cite{OG3}. For the remaining cases, we have the following answer about the existence of symplectic resolution (see \cite{LS}, \cite{KLS}):

\begin{thm}
\label{thm:ls}Let $S$ be an abelian or projective K3 surface, $v=mw$ a Mukai vector such that $m\geq 2$ and $w^{2}=2k$ for some $k\in\mathbb{N}$. If $w=(r,\xi,a)$, suppose that either $r> 0$ and $\xi\in NS(S)$, or $r=0$, $\xi$ is the first Chern class of an effective divisor and $a\neq 0$. Finally, let $H$ be a $v-$generic polarization. Then:
\begin{enumerate}
 \item if $m=2$ and $w^{2}=2$, then $M_{v}(S,H)$ admits a symplectic resolution $\pi_{v}:\widetilde{M}_{v}=\widetilde{M}_{v}(S,H)\longrightarrow M_{v}$, obtained as the blow-up of $M_{v}$ along the singular locus $\Sigma_{v}=M_{v}\setminus M_{v}^{s}$ with reduced structure (\textbf{Lehn-Sorger});
 \item if $m\geq 3$, or $m=2$ and $k\geq 2$, then $M_{v}(S,H)$ does not admit any symplectic resolution and it is locally factorial (\textbf{Kaledin-Lehn-Sorger}).
\end{enumerate}
\end{thm}

In this paper, we deal with the moduli spaces verifying the conditions of point 1 of Theorem \ref{thm:ls}. We resume these conditions in the following:

\begin{defn}
Let $S$ be an abelian or projective K3 surface, $v$ a Mukai vector, $H$ an ample line bundle on $S$. We say that $(S,v,H)$ is an \textit{OLS-triple} if the following conditions are verified:
\begin{enumerate}
 \item the polarization $H$ is primitive and $v-$generic as in Definition \ref{defn:our};
 \item we have $v=2w$, where $w$ is a primitive Mukai vector such that $w^{2}=2$;
 \item if $w=(r,\xi,a)$, we have $r\geq 0$, $\xi\in NS(S)$, and if $r=0$ then $\xi$ is the first Chern class of an effective divisor.
\end{enumerate}
\end{defn}

The name OLS-triple is chosen because they were first studied by O'Grady in \cite{OG2}, \cite{OG3} and Lehn-Sorger in \cite{LS}. If $(S,v,H)$ is an OLS-triple, then $M_{v}(S,H)$ admits a symplectic resolution $\widetilde{M}_{v}(S,H)$ obtained by blowing-up $M_{v}(S,H)$ along its singular locus with reduced structure. If $S$ is abelian, let $$\widetilde{K}_{v}=\widetilde{K}_{v}(S,H):=\pi_{v}^{-1}(K_{v}),$$and we still write $\pi_{v}:\widetilde{K}_{v}\longrightarrow K_{v}$ for the symplectic resolution. 

The aim of the present paper it to generalize Theorem \ref{thm:vprimitivo} to OLS-triples. Namely, we prove the two following theorems, which are the main results of this paper:

\begin{thm}
\label{thm:main}Let $(S,v,H)$ be an OLS-triple.
\begin{enumerate}
 \item If $S$ is K3, then $\widetilde{M}_{v}(S,H)$ is an irreducible symplectic variety which is deformation equivalent to $\widetilde{M}_{10}$.
 \item If $S$ is abelian, then $\widetilde{K}_{v}(S,H)$ is an irreducible symplectic variety which is deformation equivalent to $\widetilde{K}_{6}$.
\end{enumerate}
\end{thm}

The proof of this Theorem is contained in Section 2. The idea is to use deformations of the moduli spaces and of their symplectic resolutions induced by deformations of the underlying surfaces, and isomorphisms between moduli spaces with different Mukai vectors which are induced by Fourier-Mukai transforms (the main ingredient here is given by some results of Yoshioka, see \cite{Y3}).

\begin{thm}
\label{thm:mainl}Let $(S,v,H)$ be an OLS-triple. 
\begin{enumerate}
 \item If $S$ is K3, then $\pi_{v}^{*}:H^{2}(M_{v},\mathbb{Z})\longrightarrow H^{2}(\widetilde{M}_{v},\mathbb{Z})$ is injective, and the restrictions to $H^{2}(M_{v},\mathbb{Z})$ of the pure weight-two Hodge structure and of the Beauville form on $H^{2}(\widetilde{M}_{v},\mathbb{Z})$ give a pure weight-two Hodge structure and a lattice structure on $H^{2}(M_{v},\mathbb{Z})$. Moreover, there is a Hodge isometry $$\lambda_{v}:v^{\perp}\longrightarrow H^{2}(M_{v},\mathbb{Z}).$$
 \item If $S$ is abelian, then $\pi_{v}^{*}:H^{2}(K_{v},\mathbb{Z})\longrightarrow H^{2}(\widetilde{K}_{v},\mathbb{Z})$ is injective, and the restrictions to $H^{2}(K_{v},\mathbb{Z})$ of the pure weight-two Hodge structure and of the Beauville form on $H^{2}(\widetilde{K}_{v},\mathbb{Z})$ give a pure weight-two Hodge structure and a lattice structure on $H^{2}(K_{v},\mathbb{Z})$. Moreover, there is a Hodge isometry $$\nu_{v}:v^{\perp}\longrightarrow H^{2}(K_{v},\mathbb{Z}).$$
\end{enumerate}
\end{thm}

The proof of this is contained in Section 3. The reason why $\pi_{v}^{*}$ is injective is because the singularities of $M_{v}$ and $K_{v}$ are rational. The construction of the morphism $\lambda_{v}$ is a generalization of that of the Mukai-Donaldson morphism. Using Theorem \ref{thm:main} and some commutativity of diagrams we can reduce to the case of $M_{10}$ or $K_{6}$: there one finally uses results of \cite{Pe} and \cite{R1} to conclude.

\section{Deformations of moduli spaces}

In this section we study how moduli spaces and their symplectic resolutions vary under deformation. In section 2.1 we recall the notion of $v-$genericity of polarizations. In section 2.2 we introduce the main deformation we will look at, i. e. the deformation of a moduli space and of its symplectic resolution induced by a deformation of an OLS-triple along a smooth, connected curve. In section 2.3 we give explicit deformations of OLS-triples whose Mukai vector has positive rank, and in section 2.4 we use these and some results of \cite{Y3} to prove Theorem \ref{thm:main}.

\subsection{The $v-$genericity and properties}

In this section we recall the definition of $v-$generic polarization. Let $S$ be an abelian or projective K3 surface, $v=(v_{0},v_{1},v_{2})$ a Mukai vector on $S$. Moreover, let $H$ be any polarization, i. e. $H\in Amp(S)$, where $Amp(S)\subseteq NS(S)$ is the ample cone of $S$.

\begin{defn}
\label{defn:our}
A polarization $H$ is $v-$\textit{generic} if for every $H-$polystable sheaf $\mathscr{E}$ of Mukai vector $v$ and every direct summand $\mathscr{F}$ of $\mathscr{E}$, we have $v(\mathscr{F})\in\mathbb{Q}\cdot v$.
\end{defn}

We give a characterization of the $v-$genericity of a polarization. Let $H$ be any ample divisor on $S$, and let $v=(v_{0},v_{1},v_{2})$ be a Mukai vector. Let $\mathscr{E}$ be an $H-$semistable sheaf with Mukai vector $v$, and let $\mathscr{F}\subseteq\mathscr{E}$ an $H-$destabilizing subsheaf with Mukai vector $u:=(u_{0},u_{1},u_{2})$.

\begin{defn}
The \textit{divisor associated to the pair} $(\mathscr{E},\mathscr{F})$ is defined as follows:
\begin{enumerate}
 \item if $v_{0}>0$, it is the divisor $D:=u_{0}v_{1}-v_{0}u_{1}$;
 \item if $v_{0}=0$, it is the divisor $D:=u_{2}v_{1}-v_{2}u_{1}$.
\end{enumerate}
The set of the non-zero divisors associated to all the possible pairs is denoted $W_{v}(H)$.
\end{defn}

The characterisation of the $v-$genericity is the following:

\begin{lem}
\label{lem:genpoly}Let $v=(v_{0},v_{1},v_{2})$ be a Mukai vector such that if $v_{0}=0$, then $v_{2}\neq 0$. A polarization $H$ is $v-$generic if and only if $W_{v}(H)=\emptyset$.
\end{lem}

\proof We present a proof for $v_{0}>0$, as the one for $v_{0}=0$ is similar. Suppose that $W_{v}(H)=\emptyset$, and let $\mathscr{E}$ be a strictly $H-$polystable sheaf with Mukai vector $v$. Let $\mathscr{F}$ be a direct summand of $\mathscr{E}$, and let $v(\mathscr{F})=(u_{0},u_{1},u_{2})$. Recall that $$p_{H}(\mathscr{E})=n^{2}+2\frac{v_{1}\cdot H}{v_{0}H^{2}}n+2\frac{v_{0}+v_{2}}{v_{0}H^{2}},\,\,\,\,\,\,\,p_{H}(\mathscr{F})=n^{2}+2\frac{u_{1}\cdot H}{u_{0}H^{2}}n+2\frac{u_{0}+u_{2}}{u_{0}H^{2}},$$and as $\mathscr{F}$ is a direct summand of $\mathscr{E}$, we have $p_{H}(\mathscr{E})=p_{H}(\mathscr{F})$. This gives $u_{2}=(u_{0}/v_{0})v_{2}$ and as $W_{v}(H)=\emptyset$, the divisor $D$ associated to the pair $(\mathscr{E},\mathscr{F})$ is trivial, so that $u_{1}=(u_{0}/v_{0})v_{1}$. Hence $u=(u_{0}/v_{0})v$, and $H$ is $v-$generic.

For the converse, suppose that $H$ is $v-$generic. Let $\mathscr{E}$ be an $H-$semistable sheaf with Mukai vector $v$, and let $\mathscr{F}$ be any $H-$destabilizing subsheaf of $\mathscr{E}$ with Mukai vector $u=(u_{0},u_{1},u_{2})$. The $v-$genericity of $H$ implies that $u=(p/q)v$ for some $p,q\in\mathbb{Z}$. This means that $p/q=u_{0}/v_{0}$, and that $u_{1}=(u_{0}/v_{0})v_{1}$. Then the divisor associated to the pair $(\mathscr{E},\mathscr{F})$ is trivial, and we are done. \endproof

We now study the relation between the notion of $v-$genericity given in Definition \ref{defn:our} and those of \cite{Y2} and of Appendix 4.C of \cite{HL}. To do so, we recall the notions of $v-$walls and $v-$chambers: as these notions when $v_{0}>0$ are different than when $v_{0}=0$, we need to divide in these two cases.

\subsubsection{Walls and chambers for $v_{0}\geq 2$}

In this section we suppose $v_{0}\geq 2$. If $S$ is a K3 surface, let $$|v|:=\frac{v^{2}_{0}}{4}(v,v)+\frac{v^{4}_{0}}{2}.$$If $S$ is abelian, let $$|v|:=\frac{v^{2}_{0}}{4}(v,v)+\frac{v^{2}_{0}}{2}.$$Notice that $|v|$ depends only on $(v,v)$ and $v_{0}$, and as $v_{0}\geq 2$, then $|v|>0$. Hence it makes sense to define $$W_{v}:=\{D\in NS(S)\,|\,-|v|\leq D^{2}<0\}.$$By Theorem 4.C.3 of \cite{HL}, we have $W_{v}(H)\subseteq W_{v}$ for every $H\in Amp(S)$.

\begin{defn}
Let $D\in W_{v}$. The $v-$\textit{wall associated to} $D$ is $$W^{D}:=\{\alpha\in Amp(S)\,|\,D\cdot\alpha=0\}.$$
\end{defn}

Notice that the $v-$wall associated to $D\in W_{v}$ is a hyperplane in $Amp(S)$. By Theorem 4.C.2 of \cite{HL} the subset $\bigcup_{D\in W_{v}}W^{D}\subseteq Amp(S)$ is locally finite.

\begin{defn}
A connected component of $Amp(S)\setminus\bigcup_{D\in W_{v}}W^{D}$ is called $v-$\textit{chamber}.
\end{defn}

A $v-$chamber is then an open connected subcone of $Amp(S)$. Now, these $v-$chambers are important as if a polarization is in a $v-$chamber, then it is $v-$generic, as shown in the following:

\begin{lem}
\label{lem:chambervgen}Let $v=(v_{0},v_{1},v_{2})$ be a Mukai vector such that $v_{0}\geq 2$, and let $\mathcal{C}$ be a $v-$chamber. If $H\in\mathcal{C}$, then $H$ is $v-$generic. 
\end{lem}

\proof If $H\in\mathcal{C}$, then $H\notin W^{D}$, i. e. $D\cdot H\neq 0$ for every $D\in W_{v}$. As $W_{v}(H)\subseteq W_{v}$, we have that $H$ is $v-$generic.\endproof

If $H$ is a $v-$generic polarization, then it is not necessarily contained in some $v-$chamber. It is easy to produce examples of $v-$generic polarizations lying on some $v-$wall, as the following:

\begin{esem}
\label{esem:hsumuro}{\rm Let $Y$ be an elliptic K3 surface with $NS(Y)=\mathbb{Z}\cdot\sigma\oplus\mathbb{Z}\cdot f$, where $\sigma$ is a section and $f$ is a fiber. The line bundle $H=\sigma+3f$ is ample, hence it gives a polarization on $Y$. Consider the Mukai vector $v:=(2,\sigma+2f,a)$, where $a$ is odd and $a\leq 1$. Then $v$ is primitive, $v^{2}\geq -2$ and $|v|\geq 6$. Notice that $D:=\sigma-f$ is such that $D^{2}=-4$, so that $D\in W_{v}$, and $D\cdot H=0$, so that $H$ is not in a $v-$chamber. Nevertheless, $H$ is $v-$generic: if $\mathscr{E}$ is an $H-$semistable sheaf of Mukai vector $v$, then $\chi(\mathscr{E})=a+2$ is odd, so that $\mathscr{E}$ cannot be strictly $H-$semistable: if it was not $H-$stable, then $\chi(\mathscr{E})$ would be even. Hence $W_{v}(H)=\emptyset$, so that $H$ is $v-$generic by Lemma \ref{lem:genpoly}.}
\end{esem}

In general, the moduli space $M_{v}(S,H)$ depends on the choice of $H$. We now show that $M_{v}(S,H)$ does not change when $H$ is a $v-$generic polarization moving in the closure of a $v-$chamber.

\begin{prop}
\label{prop:chamber}
Let $v=(v_{0},v_{1},v_{2})$ be a Mukai vector such that $v_{0}\geq 2$. Let $\mathcal{C}$ be a $v-$chamber, and suppose $H,H'\in\overline{\mathcal{C}}$ are two $v-$generic polarizations. Then a sheaf $\mathscr{E}$ of Mukai vector $v$ is $H-$semistable if and only if it is $H'-$semistable, i. e. there is an identification $M_{v}(S,H)=M_{v}(S,H')$.
\end{prop}

\proof We present a proof for $H'\in\mathcal{C}$ and $H\in\partial\overline{\mathcal{C}}$. The same proof works for $H,H'\in\mathcal{C}$ (see even \cite{Z} for a complete proof in this case), so that we are done. Let $$[H,H']:=\{H_{t}:=tH+(1-t)H'\,|\,t\in[0,1]\}.$$As $H'\in\mathcal{C}$ and $\mathcal{C}$ is a convex cone in $Amp(S)$, then $[H,H']\cap\partial(\overline{\mathcal{C}})=\{H\}$. We show that if a sheaf $\mathscr{E}$ is $H-$semistable and $v(\mathscr{E})=v$, then $\mathscr{E}$ is $H'-$semistable. The statement then follows: we have an inclusion $f:M_{v}(S,H)\longrightarrow M_{v}(S,H')$. As $H'\in\mathcal{C}$, then $M_{v}(S,H')$ and $M_{v}(S,H)$ are both irreducible, reduced, normal and have the same dimension, hence $f$ is an isomorphism. Recall that if $\mathscr{G}$ is a sheaf with $v(\mathscr{G})=(v_{0},v_{1},v_{2})$, $v_{0}>0$ and $L$ is a polarization, then $$p_{L}(\mathscr{G})=n^{2}+2\frac{v_{1}\cdot L}{v_{0}L^{2}}+2\frac{v_{2}+v_{0}}{v_{0}L^{2}}.$$

Let $\mathscr{E}$ be $H-$semistable with $v(\mathscr{E})=(v_{0},v_{1},v_{2})$, $v_{0}>0$, and suppose that $\mathscr{E}$ is not $H'-$semistable. Then there is $\mathscr{F}\subseteq\mathscr{E}$ with $v(\mathscr{F}):=(u_{0},u_{1},u_{2})$ such that $p_{H'}(\mathscr{F})>p_{H'}(\mathscr{E})$, hence $u_{0}>0$. Suppose that there are $c,d\in\mathbb{Z}$ such that $cu_{1}=dv_{1}$. For every polarization $L$ we then have
\begin{equation}
\label{eq:cc}
\frac{u_{1}\cdot L}{u_{0}}=\bigg(\frac{dv_{0}}{cu_{0}}\bigg)\frac{v_{1}\cdot L}{v_{0}}.
\end{equation}
As $\mathscr{F}$ is $H'-$destabilizing, we have either $u_{1}\cdot H'/u_{0}>v_{1}\cdot H'/v_{0}$, or $u_{1}\cdot H'/u_{0}=v_{1}\cdot H'/v_{0}$ and $u_{2}/u_{0}>v_{2}/v_{0}$. The first condition implies $dv_{0}/cu_{0}>1$. Hence $u_{1}\cdot H/u_{0}>v_{1}\cdot H/v_{0}$ by equation (\ref{eq:cc}), which is not possible as $\mathscr{E}$ is $H-$semistable. The second condition implies $du_{0}/cv_{0}=1$, so that $u_{1}\cdot H/u_{0}=v_{1}\cdot H/v_{0}$. As $u_{2}/u_{0}>v_{2}/v_{0}$, we get $p_{H}(\mathscr{F})>p_{H}(\mathscr{E})$, which again is not possible.

In conclusion, $u_{1}$ is not a rational multiple of $v_{1}$. As $\mathscr{E}$ is $H-$semistable, we have $u_{1}\cdot H/u_{0}\leq v_{1}\cdot H/v_{0}$. Let $D$ be the divisor associated to the pair $(\mathscr{E},\mathscr{F})$. If $u_{1}\cdot H/u_{0}=v_{1}\cdot H/v_{0}$, then we have $D\cdot H=0$. As $H$ is $v-$generic, we have $D=0$, i. e. $u_{1}=(u_{0}/v_{0})v_{1}$, which is not possible as $u_{1}$ is not a rational multiple of $v_{1}$. In conclusion, we need $u_{1}\cdot H/u_{0}<v_{1}\cdot H/v_{0}$. As $\mathscr{F}$ is $H'-$destabilizing, then $u_{1}\cdot H'/u_{0}\geq v_{1}\cdot H'/v_{0}$. The function $$f:[0,1]\longrightarrow\mathbb{R},\,\,\,\,\,\,\,\,\,\,\,\,f(t):=D\cdot H_{t}$$is continuous with $f(0)<0$, $f(1)\geq 0$. Hence there is $t\in(0,1]$ such that $f(t)=0$. As $D\cdot H\neq 0$, we get $D\neq 0$, so that $D\in W_{v}$. Moreover, we have $D\cdot H_{t}=0$, so that $H_{t}\in W^{D}$. As $t>0$, we have then $H_{t}\in\mathcal{C}$, hence we get a contradiction, and we are done.\endproof

\begin{oss}
\label{oss:sr1}{\rm If $(S,v,H)$ is an OLS-triple where $v=2(r,\xi,a)$ with $r>0$, then $M_{v}$ is non-empty (this was shown by Yoshioka, as signaled in \cite{KLS}) and it is a $10-$dimensional irreducible, normal, projective variety whose singular locus is $\Sigma_{v}:=M_{v}\setminus M_{v}^{s}$, and there is a symplectic resolution $\pi_{v}:\widetilde{M}_{v}\longrightarrow M_{v}$ obtained as blow-up of $M_{v}$ along $\Sigma_{v}$ with reduced structure. This is contained in \cite{LS}: even if formally stated by asking that $H$ is in a $v-$chamber, the argument applies to every $v-$generic polarization. Moreover, by Proposition \ref{prop:chamber} we can always assume that $H$ is in a $v-$chamber.}
\end{oss}

\subsubsection{Walls and chambers for $v_{0}=0$}

In this section we deal with the case $v=(0,v_{1},v_{2})$. The definition of $v-$generic polarization we gave is, for rank 0, formally the same as the one for positive rank. In any case, as Yoshioka first remarked, if $v_{0}=0$ it can happen that there is no $v-$generic polarization at all, as shown in the following:

\begin{esem}
{\rm Let $v:=(0,c_{1}(L),0)$ be primitive, and suppose that there are two distinct irreducible curves $C_{1},C_{2}\in |L|$. For $i=1,2$ let $L_{i}\in Pic(C_{i})$ be a line bundle of degree $d_{i}=g(C_{i})-1$, and let $j_{i}$ be the inclusion of $C_{i}$ in $S$. Then $j_{i*}L_{i}$ has Mukai vector $v_{i}=(0,C_{i},0)$, and it is $H-$stable for every $H\in Amp(S)$. Moreover $j_{1*}L_{1}\oplus j_{2*}L_{2}$ is strictly $H-$semistable of Mukai vector $v$, and $H$ is not $v-$generic.}
\end{esem}

We now define walls and chambers for $v=(0,v_{1},v_{2})$ and $v_{2}\neq 0$.

\begin{defn}
Let $\mathscr{E}$ be any pure sheaf with Mukai vector $v$, and let $\mathscr{F}\subseteq\mathscr{E}$ with Mukai vector $u=(0,u_{1},u_{2})$. The \textit{divisor associated to the pair} $(\mathscr{E},\mathscr{F})$ is $D:=u_{2}v_{1}-v_{2}u_{1}$. The set of the non-zero divisors associated to all the possible pairs is denoted $W_{v}$. 
\end{defn}

\begin{defn}
Let $D\in W_{v}$. The $v-$\textit{wall associated to} $D$ is $$W^{D}:=\{\alpha\in Amp(S)\,|\,\alpha\cdot D=0\}.$$
\end{defn}

As shown in \cite{Y2}, the set of $v-$walls is finite.

\begin{defn}
A connected component of $Amp(S)\setminus\bigcup_{D\in W_{v}}W^{D}$ is called $v-$\textit{chamber}.
\end{defn}

Again, if a polarization is in a $v-$chamber, then it is $v-$generic:

\begin{lem}
\label{lem:genpoly0}Let $v=(0,v_{1},v_{2})$ be a Mukai vector such that $v_{2}\neq 0$, and let $\mathcal{C}$ be a $v-$chamber. If $H\in\mathcal{C}$, then $H$ is $v-$generic. 
\end{lem}

\proof The proof is the same as the one of Lemma \ref{lem:chambervgen}.\endproof

As in the previous case, we now show that the moduli space $M_{v}(S,H)$ does not change when $H$ is a $v-$generic polarization moving inside the closure of a $v-$chamber.

\begin{prop}
\label{prop:chamber0}
Let $v=(0,v_{1},v_{2})$ be a Mukai vector such that $v_{2}\neq 0$. Let $\mathcal{C}$ be a $v-$chamber, and suppose that $H,H'\in\overline{\mathcal{C}}$ are two $v-$generic polarizations. Then a sheaf $\mathscr{E}$ of Mukai vector $v$ is $H-$semistable if and only if it is $H'-$semistable, i. e. there is an identification $M_{v}(S,H)=M_{v}(S,H')$.
\end{prop}

\proof The proof follows the same pattern as the one of Proposition \ref{prop:chamber}.\endproof

\begin{oss}
{\rm Let $(S,v,H)$ be an OLS-triple where $v=2(0,\xi,a)$, then $M_{v}$ is non-empty (this was shown by Yoshioka, as signaled in \cite{KLS}), and it is a $10-$dimensional irreducible, normal, projective variety whose singular locus is $\Sigma_{v}:=M_{v}\setminus M_{v}^{s}$, and there is a symplectic resolution $\pi_{v}:\widetilde{M}_{v}\longrightarrow M_{v}$ obtained as blow-up of $M_{v}$ along $\Sigma_{v}$ with reduced structure (see Remark \ref{oss:sr1}). Moreover, if $v=2(0,\xi,0)$ the tensorization with $H$ gives an isomorphism between $M_{v}(S,H)$ and $M_{v'}(S,H)$, where $v'=2(0,\xi,\xi\cdot H)$, where $\xi\cdot H>0$ (see Theorem 3.14 of \cite{Z}): hence, we can always suppose $v=2(0,\xi,a)$ for $a\neq 0$, and by Proposition \ref{prop:chamber0} that $H$ is in a $v-$chamber.}
\end{oss}

\subsection{Deformations of OLS-triples}

We introduce the main construction we use in the following. Let $(S,v,H)$ be an OLS-triple and $T$ a smooth, connected curve, and use the following notation: if $f:Y\longrightarrow T$ is a morphism and $\mathscr{L}\in Pic(Y)$, for every $t\in T$ we let $Y_{t}:=f^{-1}(t)$ and $\mathscr{L}_{t}:=\mathscr{L}_{|Y_{t}}$.

\begin{defn}
Let $(S,v,H)$ be an OLS-triple, where $v=2(r,\xi,a)$ and $\xi=c_{1}(L)$. A \textit{deformation} of $(S,v,H)$ on $T$ is a triple $(\mathscr{X},\mathscr{H},\mathscr{L})$, where:
\begin{enumerate}
\item $\mathscr{X}$ is a projective, smooth deformation of $S$ along $T$, i. e. there is a smooth, projective, surjective map $f:\mathscr{X}\longrightarrow T$ such that $\mathscr{X}_{t}$ is a projective K3 surface (resp. an abelian surface) for every $t\in T$, and there is $0\in T$ such that $\mathscr{X}_{0}\simeq S$;
\item $\mathscr{H}$ is a line bundle on $\mathscr{X}$ such that $\mathscr{H}_{t}$ is ample for every $t\in T$ and such that $\mathscr{H}_{0}\simeq H$;
\item $\mathscr{L}$ is a line bundle on $\mathscr{X}$ such that $\mathscr{L}_{0}\simeq L$.
\item if for every $t\in T$ we let $\xi_{t}:=c_{1}(\mathscr{L}_{t})$, $w_{t}:=(r,\xi_{t},a)$ and $v_{t}:=2w_{t}$, then we ask that $\mathscr{H}_{t}$ is $v_{t}-$generic.
\end{enumerate}
\end{defn}

\begin{oss}
\label{oss:vgenv}{\rm Notice that if $(S,v,H)$ is an OLS-triple and $(\mathscr{X},\mathscr{L},\mathscr{H})$ is a deformation of $(S,v,H)$ along a smooth, connected curve $T$, then $(\mathscr{X}_{t},v_{t},\mathscr{H}_{t})$ is an OLS-triple for every $t\in T$. Indeed, we have $v_{t}=2w_{t}$, where $w_{t}$ is primitive and $w^{2}_{t}=2$. Moreover, if $r=0$, then $\xi_{t}$ is effective: we have $\xi_{t}^{2}=2$, hence either $\xi_{t}$ or $-\xi_{t}$ is effective; as $\xi$ is effective, then $-\xi\cdot H<0$, so that $-\xi_{t}\cdot\mathscr{H}_{t}<0$, hence $\xi_{t}$ is effective.}
\end{oss}

\begin{oss}
\label{oss:ap}{\rm Consider an OLS-triple $(S,v,H)$ where $S$ is a projective K3 surface (resp. an abelian surface) and $v=2(r,\xi,a)$, where $\xi=c_{1}(L)$. Let $T$ be a smooth, connected curve. Moreover, consider a smooth, projective deformation $f:\mathscr{X}\longrightarrow T$ of $S$ such that $\mathscr{X}_{0}\simeq S$, and on $\mathscr{X}$ consider two line bundles $\mathscr{H}$ and $\mathscr{L}$ such that $\mathscr{H}_{0}\simeq H$ and $\mathscr{L}_{0}\simeq L$. In general, we have that $(\mathscr{X}_{t},v_{t},\mathscr{H}_{t})$ is not an OLS-triple for every $t\in T$: by Remark \ref{oss:vgenv}, this is the case if and only if $\mathscr{H}_{t}$ is ample and $v_{t}-$generic for every $t\in T$. Up to removing a finite number of points from $T$, we can always assume that this is the case. Indeed, the set of $t\in T$ such that $\mathscr{H}_{t}$ is not ample is finite. For the $v_{t}-$genericity, we have the following Proposition.}
\end{oss}

\begin{prop}
\label{prop:apertura}Let $(S,v,H)$ be an OLS-triple where $S$ is a projective K3 surface (resp. an abelian surface) and $v=2(r,\xi,a)$, $\xi=c_{1}(L)$. Let $T$ be a smooth, connected curve, $f:\mathscr{X}\longrightarrow T$ a smooth, projective deformation of $S$ such that $\mathscr{X}_{0}\simeq S$, $\mathscr{H}\in Pic(\mathscr{X})$ such that $\mathscr{H}_{t}$ is ample for every $t\in T$ and $\mathscr{H}_{0}\simeq H$, and $\mathscr{L}\in Pic(\mathscr{X})$ such that $\mathscr{L}_{0}\simeq L$. Then the set $$T':=\{t\in T\,|\,\mathscr{H}_{t}\,\,is\,\,not\,\,v_{t}-generic\}$$is finite.
\end{prop}

\proof Consider the relative moduli space of semistable (resp. stable) sheaves $$\phi:\mathscr{M}\longrightarrow T,\,\,\,\,\,\,\,\,\,\,\,\,\,\,\,(resp.\,\,\,\,\,\,\,\phi^{s}:\mathscr{M}^{s}\longrightarrow T)$$associated to $f$, where $\phi$ (resp. $\phi^{s}$) is a projective (resp. quasi-projective) map, such that $\mathscr{M}_{t}=M_{v_{t}}(\mathscr{X}_{t},\mathscr{H}_{t})$ (resp. $\mathscr{M}_{t}^{s}=M_{v_{t}}^{s}(\mathscr{X}_{t},\mathscr{H}_{t})$) for every $t\in T$. By Theorem 4.3.7 in \cite{HL}, $\mathscr{M}^{s}$ is an open subset of $\mathscr{M}$. Let $\mathscr{M}^{ss}:=\mathscr{M}\setminus\mathscr{M}^{s}$, which is a closed subset of $\mathscr{M}$. By definition, $\mathscr{M}^{ss}_{t}$ parameterizes strictly $\mathscr{H}_{t}-$polystable sheaves for every $t\in T$. There is an irreducible component $\Sigma$ of $\mathscr{M}^{ss}$ such that $\Sigma\cap\mathscr{M}_{t}$ parameterizes strictly $\mathscr{H}_{t}-$polystable sheaves whose direct summands are two $\mathscr{H}_{t}-$stable sheaves with Mukai vector $w_{t}$. Let $\Xi$ be the union of all the other irreducible components of $\mathscr{M}^{ss}$. By definition of $v-$genericity we have $T'=\phi(\Xi)$. As $\Xi$ is closed in $\mathscr{M}$ and the morphism $\phi$ is projective, then $\phi(\Xi)$ is a closed subset of $T$. Moreover, notice that as $(\mathscr{X}_{0},v_{0},\mathscr{H}_{0})=(S,v,H)$ is an OLS-triple, then $\mathscr{H}_{0}$ is $v_{0}-$generic. Then $0\notin T'$, and $T'$ is a proper closed subset of $T$.\endproof

The reason why we introduce the notion of deformation of an OLS-triple, is because they allow us to study how the algebraic structure of the corresponding moduli space (and of its symplectic resolution) varies under variations of the algebraic structure of the base surface. The first result we prove is that the relative moduli space $\phi:\mathscr{M}\longrightarrow T$ associated to a deformation of an OLS-triple along a smooth, connected curve $T$, is a flat family.

\begin{lem}
\label{lem:redflat}
Let $(S,v,H)$ be an OLS-triple, $T$ a smooth, connected curve, and $(\mathscr{X},\mathscr{L},\mathscr{H})$ a deformation of $(S,v,H)$ along $T$. Then $\phi:\mathscr{M}\longrightarrow T$ is flat. 
\end{lem}

\proof Let $t\in T$, $T^{0}:=T\setminus\{t\}$ and $\mathscr{M}^{0}:=\phi^{-1}(T^{0})$. The morphism $\phi$ is flat over $t$ if and only if the fiber $\mathscr{M}_{t}$ is the limit of the fibers $\mathscr{M}_{s}$ as $s\rightarrow t$, by Lemma II-29 of \cite{EH}. Now, the limit is the fiber over $t$ of the closure of the family $\mathscr{M}^{0}$, hence there is an inclusion of the limit in $\mathscr{M}_{t}$. Recall that $\mathscr{M}_{t}=M_{v_{t}}(\mathscr{X}_{t},\mathscr{H}_{t})$ is reduced and irreducible, hence it has to coincide with the previous limit.\endproof

If $(S,v,H)$ is an OLS-triple, then choosing a non-trivial deformation of it along a smooth, connected curve $T$ we get a flat, projective deformation $\phi:\mathscr{M}\longrightarrow T$ of $M_{v}(S,H)$. We now present the main result of this section, which is about local properties of this deformation: it is easy to see that if $t_{0}\in T$ is any point and $U$ is any open neighborhood of $t_{0}$ in $T$, then $\phi^{-1}(U)$ is not isomorphic to a product $\mathscr{M}_{t_{0}}\times U$. However, we show in the following Proposition, that this is true locally at every point in every fibre. 

\begin{prop}
\label{prop:loctrivial}Let $(S,v,H)$ be an OLS-triple, $T$ a smooth, connected curve, and $(\mathscr{X},\mathscr{L},\mathscr{H})$ a deformation of $(S,v,H)$ along $T$. Let $p\in\mathscr{M}$ and $0:=\phi(p)$. Then the germ $(\mathscr{M},p)$ is isomorphic, as germ of complex spaces, to the product $(T,0)\times(\mathscr{M}_{0},p)$.
\end{prop}

\proof We need the following definition: let $\phi:\mathscr{M}\longrightarrow T$ (resp. $\phi^{s}:\mathscr{M}^{s}\longrightarrow T$) be the relative moduli space of semistable (resp. stable) sheaves associated to the deformation $(\mathscr{X},\mathscr{L},\mathscr{H})$ of $(S,v,H)$ along $T$. Let $\Sigma:=\mathscr{M}\setminus\mathscr{M}^{s}$, which is a closed subset of $\mathscr{M}$. Notice that $$\Sigma=\bigcup_{t\in T}\Sigma_{v_{t}},$$and we use the notation $\Sigma_{t}:=\Sigma\cap\mathscr{M}=\Sigma_{v_{t}}$. MOreover, for every $t\in T$ let $\Omega_{t}$ be the singular locus of $\Sigma_{t}$, and $\Omega$ be the closed subset of $\mathscr{M}$ parameterizing sheaves of the form $\mathscr{E}\oplus\mathscr{E}$, where $\mathscr{E}$ is stable. Notice that $$\Omega=\bigcup_{t\in T}\Omega_{t}.$$As $\mathscr{M}_{t}=M_{v_{t}}(\mathscr{X}_{t},\mathscr{H}_{t})$ for every $t\in T$, the point $p\in\mathscr{M}_{0}$ is one of the following: $p$ is smooth, i. e. $p\in\mathscr{M}^{s}_{0}$; $p\in\Sigma_{0}\setminus\Omega_{0}$, i. e. $p$ is a singular point of type $A_{1}$; $p\in\Omega_{0}$. If $p$ is smooth, the result is trivial, and there is nothing to prove. Hence, we suppose $p\in\Sigma_{0}$. We have then two cases:

\textit{Case 1:} $p\in\Sigma_{0}\setminus\Omega_{0}$. Consider the Zariski tangent space $T_{p}\mathscr{M}$ of $\mathscr{M}$ at $p$: we have $dim(T_{p}\mathscr{M})=dim(\mathscr{M})+1=12$. Consider an analytic open neighborhood $U$ of $p$ in $T_{p}\mathscr{M}$, so we can view it as an open neighborhood of $0$ in $\mathbb{C}^{12}$. We let $x_{1},...,x_{12}$ be a coordinate system on $\mathbb{C}^{12}$: we can suppose $x_{12}=t$, a local coordinate of $T$ at $0$, and the point $p$ corresponds to the point $(0,...,0)$. Moreover, $U\cap\mathscr{M}$ is an analytic subvariety of $U$ of codimension 1, hence there is $f\in\mathscr{O}^{hol}(U)$, a holomorphic function on $U$, such that the equation of $U\cap\mathscr{M}$ is $f(x_{1},...,x_{11},t)=0$. Finally, we can choose $U$ so that $U\cap\Omega=\emptyset$. By the functoriality properties of the relative moduli space, we have that $\Sigma\setminus\Omega$ is smooth and submersive on $T$, so that we can suppose that the equation of $U\cap\Sigma$ is $x_{1}=x_{2}=x_{3}=0$. As near $p$ the fibre $\mathscr{M}_{0}$ is analytically isomorphic to a product of an $A_{1}-$singularity with a smooth polydisc (see \cite{LS}), we have $$f(x_{1},...,x_{11},0)=x_{1}^{2}+x_{2}^{2}+x_{3}^{2},$$so that $$f(x_{1},...,x_{11},t)=x_{1}^{2}+x_{2}^{2}+x_{3}^{2}+\sum_{j}p_{j}(x_{1},...,x_{11})t^{j},$$where $p_{j}$ are polynomials on $U$ depending only on $x_{1},...,x_{11}$. Moreover, we have $p_{j}\in I^{2}$ for every $j$, where $I$ is the ideal of $\mathscr{O}^{hol}(U)$ generated by $x_{1},x_{2}$ and $x_{3}$.
 
Now, let $p:U\longrightarrow T$ defined as $p(x_{1},...,x_{11},t):=t$, and let $V:=U\cap\Sigma$, on which we have coordinates $x_{4},...,x_{11},t$. Finally, let $$p':U\longrightarrow V,\,\,\,\,\,\,\,\,\,\,\,\,\,\,p'(x_{1},...,x_{11},t):=(x_{4},...,x_{11},t)$$and $q:V\longrightarrow T$ be defined as $q(x_{4},...,x_{11},t):=t$. Notice that $q\circ p'=p$. Moreover, the fibers of $p'$ are all singular, and the fiber over 0 has an $A_{1}-$singularity. By the deformation theory of $A_{1}-$singularities of \cite{KM}, we have that there is an open neighborhood $U'\subseteq\mathscr{M}$ of the point $p$ which is a product of an $A_{1}-$singularity by a $9-$dimensional polydisc $D$. As $\phi_{|U'}$ is identified with $p_{|U'}$, then the projection onto $D$ factors $\phi$. Hence $\phi$ is locally trivial at $p$, and we are done.

\textit{Case 2:} $p\in\Omega_{0}$. The strategy is the following: first, we show that for every $n\in\mathbb{N}$, the infinitesimal $n-$th order deformation of $\mathscr{M}_{0}$ induced by $\mathscr{M}$ is locally trivial at $p$. Once this is shown, the statement follows in this way: by Corollary 0.2 of \cite{FK} there is a maximal subspace $(T',0)\subseteq (T,0)$ such that $(\mathscr{M}_{T'},p)$ is isomorphic, as germ of complex spaces, to the product $(T',0)\times(\mathscr{M}_{0},p)$ (where $\mathscr{M}_{T'}:=\mathscr{M}\times_{T}T'$). Notice that as the infinitesimal $n-$th order deformation of $\mathscr{M}_{0}$ induced by $\mathscr{M}$ is locally trivial at $p$ for every $n$, then $T'$ is positive dimensional. As $T$ is a curve, we finally get $(T',0)=(T,0)$, and we are done.

We are left to prove that the infinitesimal $n-$th order deformation of $\mathscr{M}_{0}$ induced by $\mathscr{M}$ is locally trivial at the points of $\Omega_{0}$ for every $n$, and we proceed by induction on $n$. For $n=1$, let $$T^{1}:=\mathscr{E}xt^{1}(\Omega^{1}_{\mathscr{M}_{0}},\mathscr{O}_{\mathscr{M}_{0}}),$$where $\Omega^{1}_{\mathscr{M}_{0}}$ is the sheaf of holomorphic $1-$forms on $\mathscr{M}_{0}$: then $T^{1}$ is supported on $\Sigma_{0}$, and the local sections of $T^{1}$ correspond to local infinitesimal first order deformations of $\mathscr{M}_{0}$. Moreover, by \cite{LS} we know that $T^{1}$ is pure.

We show that the infinitesimal first order deformation of $\mathscr{M}_{0}$ induced by $\mathscr{M}$ is locally trivial at $p$: consider a Stein open neighborhood $U_{1}$ of $p$ in $\mathscr{M}_{0}$, and let $s$ be the element of $T^{1}$ on $U_{1}$ induced by $\mathscr{M}$. Let $q\in U_{1}\cap (\Sigma_{0}\setminus\Omega_{0})$: then $q$ is a point of the previous case, hence $s$ is locally trivial at $q$. This means that there is a Stein open neighborhood $V_{q}\subseteq U_{1}$ of $q$ such that $s_{|V_{q}}$ is trivial. By purity of $T^{1}$, $s$ is trivial on $U_{1}$, and we are done.

By induction, suppose that the infinitesimal $n-$th order deformation of $\mathscr{M}_{0}$ induced by $\mathscr{M}$ is locally trivial at $p$. There are two extensions of it to a local infinitesimal $(n+1)-$th order deformation at $p$: the trivial one, which we call $s_{1}$, and the one induced by $\mathscr{M}$, which we call $s_{2}$. By Theorem 2.11 of \cite{S} and Lemma 2.12 of \cite{FM} there is a transitive action of $T^{1}$ on the space of small extensions, hence there is an element $h$ of $T^{1}$ on a Stein open neighborhood $U$ of $p$ such that $h(s_{1})=s_{2}$, where $h(s_{1})$ is the action of $h$ on $s_{1}$. Let $q\in U\cap (\Sigma_{0}\setminus\Omega_{0})$: as this is a point of the previous case, the infinitesimal $(n+1)-$th order deformation of $\mathscr{M}_{0}$ induced by $\mathscr{M}$ is locally trivial at $q$, hence there is a Stein open neighborhood $V_{q}\subseteq U$ of $q$ such that $s_{2|V_{q}}=s_{1|V_{q}}$. This implies that $h_{|V_{q}}$ is trivial. Again, by purity of $T^{1}$ this implies that $h$ is trivial on $U$, so that $s_{1}=s_{2}$ on $U$, and we are done. 
\endproof

The Proposition we just proved has two important consequences. The first one is that if $(S_{1},v_{1},H_{1})$ and $(S_{2},v_{2},H_{2})$ are two OLS-triples which are related by a deformation of OLS-triples along a smooth, connected curve, then $\widetilde{M}_{v_{1}}(S_{1},H_{1})$ and $\widetilde{M}_{v_{2}}(S_{2},H_{2})$ are deformation equivalent. We now look at deformation of the symplectic resolution $\widetilde{M}_{v}(S,H)$ of the moduli space $M_{v}(S,H)$. 

\begin{prop}
\label{prop:defosvh}
Let $(S,v,H)$ be an OLS-triple, $T$ a smooth connected curve, and $(\mathscr{X},\mathscr{L},\mathscr{H})$ a deformation of $(S,v,H)$ along $T$. 
\begin{enumerate}
 \item If $S$ is a K3 surface, then $\widetilde{M}_{v}(S,H)$ is irreducible symplectic if and only if $\widetilde{M}_{v_{t}}(\mathscr{X}_{t},\mathscr{H}_{t})$ is for some (and hence for all) $t\in T$, and their deformation classes are equal.
 \item If $S$ is an abelian surface, Then $\widetilde{K}_{v}(S,H)$ is irreducible symplectic if and only if $\widetilde{K}_{v_{t}}(\mathscr{X}_{t},\mathscr{H}_{t})$ is for some (and hence for all) $t\in T$, and their deformation classes are equal.
\end{enumerate}
\end{prop}

\proof Let us suppose that $S$ is K3, and define $\pi:\widetilde{\mathscr{M}}\longrightarrow\mathscr{M}$ to be the blow-up of $\mathscr{M}$ along $\Sigma$ with reduced structure. We have a morphism $$\psi:=\phi\circ\pi:\widetilde{\mathscr{M}}\longrightarrow T,$$which is projective (as $\phi$ and $\pi$ are projective) and flat (by Lemma \ref{lem:redflat}). As an immediate consequence of Proposition \ref{prop:loctrivial}, we have that $$\widetilde{\mathscr{M}}_{t}=(Bl_{\Sigma_{red}}\mathscr{M})_{t}=Bl_{\Sigma_{t,red}}\mathscr{M}_{t}.$$Moreover, we have $\mathscr{M}_{t}=M_{v_{t}}(\mathscr{X}_{t},\mathscr{H}_{t})$ and $\Sigma_{t,red}=\Sigma_{v_{t},red}$, and $$Bl_{\Sigma_{v_{t},red}}M_{v_{t}}(\mathscr{X}_{t},\mathscr{H}_{t})=\widetilde{M}_{v_{t}}(\mathscr{X}_{t},\mathscr{H}_{t})$$for every $t\in T$. Hence $\widetilde{M}_{v_{t}}(\mathscr{X}_{t},\mathscr{H}_{t})$ is a smooth, symplectic, projective variety. As $T$ is smooth and connected, the statement follows as in the proof of Corollary 6.2.12 of \cite{HL}.

If $S$ is an abelian surface, we need one step more: define $\widehat{\mathscr{X}}:=Pic^{0}(\mathscr{X})$, with the natural map $\widehat{f}:\widehat{\mathscr{X}}\longrightarrow T$, which is again smooth. Consider $$Z:=\{(0_{\mathscr{X}_{t}},\mathscr{O}_{\mathscr{X}_{t}})\in\mathscr{X}_{t}\times\widehat{\mathscr{X}}_{t}\,|\,t\in T\}\subseteq\mathscr{X}\times_{T}\widehat{\mathscr{X}},$$with the natural morphism $g:Z\longrightarrow T$, which is clearly an isomorphism. Moreover, we can define a $T-$flat morphism $$a:\mathscr{M}\longrightarrow\mathscr{X}\times_{T}\widehat{\mathscr{X}},$$such that $a_{|\mathscr{M}_{t}}:=a_{v_{t}}$. Notice that $\phi=g\circ a$. Now, define $$\mathscr{K}:=\mathscr{M}\times_{\mathscr{X}\times_{T}\widehat{\mathscr{X}}}Z.$$Hence we have a morphism $$\phi^{0}:=g\circ p_{Z}:\mathscr{K}\longrightarrow T,$$where $p_{Z}:\mathscr{K}\longrightarrow Z$ is the natural projection. Notice that $\phi^{0}$ is flat and projective, and that for every $t\in T$ we have $\mathscr{K}_{t}=K_{v_{t}}(\mathscr{X}_{t},H_{t})$: in conclusion, $\mathscr{K}$ is a projective, flat deformation of $K_{v}$. Now, let $\widetilde{\mathscr{K}}:=\widetilde{\mathscr{M}}\times_{\mathscr{M}}\mathscr{K}$, which has a natural projection $$\psi^{0}:=\phi^{0}\circ p_{\mathscr{K}}:\widetilde{\mathscr{K}}\longrightarrow T,$$which is flat and projective. The rest of the proof now goes as in the case of K3 surfaces.\endproof

The second consequence of Proposition \ref{prop:loctrivial} is that the family $\phi:\mathscr{M}\longrightarrow T$ is topologically a product on small open subsets of $T$. More precisely, we have the following:

\begin{cor}
\label{cor:banaletop}Let $(S,v,H)$ be an OLS-triple, $T$ a smooth, connected curve and $(\mathscr{X},\mathscr{L},\mathscr{H})$ a deformation of $(S,v,H)$ along $T$. Let $\phi:\mathscr{M}\longrightarrow T$ be the relative moduli space induced by $(\mathscr{X},\mathscr{L},\mathscr{H})$. Then for every $t\in T$ there is an open neighborhood $U\subseteq T$ of $t$, and a homeomorphism $$h:\phi^{-1}(U)\longrightarrow \mathscr{M}_{t}\times U,$$such that $p_{U}\circ h=\phi$, where $p_{U}:\mathscr{M}_{t}\times U\longrightarrow U$ is the projection.
\end{cor}

\proof Let $\mathcal{S}:=\{S_{1},S_{2},S_{3}\}$ be the stratification of $\mathscr{M}$ given by $S_{1}:=\mathscr{M}\setminus\Sigma$, $S_{2}:=\Sigma\setminus\Omega$ and $S_{3}:=\Omega$. Let us suppose that the stratification $\mathcal{S}$ is Whitney: as every stratum $S_{i}$ is submersive on $T$, and as the closure $\overline{S_{i}}$ of $S_{i}$ in $\mathscr{M}$ is proper over $T$, the statement then follows from the Thom First Homotopy Lemma (see Theorem 3.5 of \cite{D}).

We then need to show that the stratification $\mathcal{S}$ is Whitney, i. e. that for every $i=1,2,3$, the stratum $S_{i}$ is Whitney regular over $S_{j}$ for every $j>i$ (see Definition 1.7 of \cite{D}). Let $t\in T$, and recall that for every $t\in T$ we have $\mathscr{M}_{t}=M_{v_{t}}$. Let $\mathcal{S}_{t}=\{S_{1,t},S_{2,t},S_{3,t}\}$ be the stratification of $M_{v_{t}}$ obtained as $S_{i,t}:=S_{i}\cap\mathscr{M}_{t}$. More explicitely, we have $S_{1,t}:=M_{v_{t}}\setminus\Sigma_{t}$, $S_{2,t}:=\Sigma_{t}\setminus\Omega_{t}$ and $S_{3,t}:=\Omega_{t}$. Notice that the stratification $\mathcal{S}$ of $\mathscr{M}$ is Whitney if the stratification $\mathcal{S}_{t}$ of $M_{v_{t}}$ is Whitney for every $t\in T$: indeed, by Proposition \ref{prop:loctrivial} we know that for every $p\in\mathscr{M}$ the germ $(\mathscr{M},p)$ is isomorphic, as germ of complex spaces, to the product $(T,\phi(p))\times(\mathscr{M}_{\phi(p)},p)$.

In conclusion, we need to show that the stratification $\mathcal{S}_{t}$ is Whitney for every $t\in T$, i. e. that for every $i=1,2,3$, the stratum $S_{i,t}$ is Whitney regular over each stratum $S_{j,t}$ with $j>i$. We have two cases:

\textit{Case 1}: $S_{1,t}$ is Whitney regular over $S_{2,t}$. To do so, let $p\in S_{2,t}=\Sigma_{t}\setminus\Omega_{t}$: then there is an open neighborhood $U\subseteq M_{v_{t}}$ of $p$, which is a product of a type $A_{1}-$singularity by an $8-$dimensional polydisc. As the stratification of the singularities of the type $A_{1}-$singularity is Whitney, this implies the Whitney regularity of $S_{1}$ over $S_{2}$. 

\textit{Case 2}: $S_{1,t}$ and $S_{2,t}$ are Whitney regular over $S_{3,t}$. Let $q\in S_{3,t}=\Omega_{t}$: then there is open neighborhood $V\subseteq M_{v_{t}}$ of $q$ which is a product of a singular $6-$dimensional variety $Z$ with a $4-$dimensional polydisc. Let $Z_{sing}$ be the singular locus of $Z$, and $(Z_{sing})_{sing}$ the singular locus of $Z_{sing}$. Notice that $Z_{sing}=\Sigma_{t}\cap V$ and $(Z_{sing})_{sing}=\Omega_{t}\cap V$. Hence, the stratification $\mathcal{S}_{t}$ induces a strtification $\mathcal{S}':=\{S'_{1},S'_{2},S'_{3}\}$ on $Z$, where $S'_{1}:=Z\setminus Z_{sing}$, $S'_{2}:=Z_{sing}\setminus(Z_{sing})_{sing}$ and $S'_{3}:=(Z_{sing})_{sing}$. To show that the stratification $\mathcal{S}_{t}$ on $M_{v_{t}}$ is Whitney, we then just need to show that the stratification $\mathcal{S}'$ of $Z$ is Whitney. Again, we show that for every $i=1,2,3$ the stratum $S'_{i}$ is Whitney regular over $S'_{j}$ for every $j>i$. The fact that $S'_{1}$ is Withney regular over $S'_{2}$ follows as in Case 1. We are then left to show that $S'_{1}$ and $S'_{2}$ are Withney regular over $S'_{3}$: but notice that as $Z$ is $6-$dimensional, we have that $(Z_{sing})_{sing}$ is $0-$dimensional, so that the stratum $S'_{3}$ is given by a point. By Lemma 1.10 of \cite{D}, we have that $S'_{1}$ and $S'_{2}$ are Whitney regular over $S'_{3}$, and we are done.\endproof  

\subsection{Deformations and Mukai vectors of positive rank}

In this section we consider Mukai vectors $v$ of positive rank, and we show that the deformation class of $\widetilde{M}_{v}$  and $\widetilde{K}_{v}$ depends only on the rank of $v$. To do so, we follow closely the arguments used by O'Grady in \cite{OG1}. As first step, we remark that the tensorization via a line bundle does not change the moduli spaces. Let $S$ be an abelian or projective K3 surface.

\begin{defn}
Let $v,v'\in\widetilde{H}(S,\mathbb{Z})$ be two Mukai vectors, $v=(v_{0},v_{1},v_{2})$, $v'=(v'_{0},v'_{1},v'_{2})$ and $v_{0},v'_{0}>0$. We say that $v$ and $v'$ are \textit{equivalent} if there is a line bundle $L$ on $S$ such that $v'=v\cdot ch(L)$.
\end{defn}

If $(S,v,H)$ and $(S,v',H)$ are two OLS-triples such that $v'=v\cdot ch(L)$ for some line bundle $L\in Pic(S)$, then the tensorization with $L$ defines an isomorphism between $M_{v}(S,H)$ and $M_{v'}(S,H)$. This is due to the following, which is Lemma 1.1 of \cite{Y2}:

\begin{lem}
\label{lem:mustbir}If $v$ is a Mukai vector of positive rank, $H$ is in a $v-$chamber, and $L\in Pic(S)$, then the tensorization with $L$ gives an isomorphism between $M_{v}(S,H)$ and $M_{v'}(S,H)$.
\end{lem}

\begin{oss}
{\rm This Lemma is originally stated only for stable sheaves, but the argument goes through for semistable sheaves. Moreover, we can suppose that $H$ is $v-$generic by Proposition \ref{prop:chamber}.}
\end{oss}

In order to give explicit deformations of an OLS-triple $(S,v,H)$ where $v=2(r,\xi,a)$ and $r>0$, we use the irreducibility of the moduli space of polarized K3 or abelian surfaces. Hence, it is useful to suppose $\xi=c_{1}(H)$, which is always possible by the following: 

\begin{lem}
\label{lem:polgenxi}Let $(S,v,H)$ be an OLS-triple where $v=2(r,\xi,a)$ is such that $r>0$. Suppose that $\rho(S)\geq 2$, and let $\mathcal{C}$ be the a $v-$chamber such that $H\in\overline{\mathcal{C}}$. Then there exists a Mukai vector $v'=2(r,\xi',a')$ such that
\begin{enumerate}
\item $v'$ is equivalent to $v$;
\item $\xi'$ is a primitive ample class lying in $\mathcal{C}$.
\end{enumerate}
Moreover, we can choose $v'$ so that $(\xi')^{2}\gg 0$.
\end{lem}

\proof First of all, by Proposition \ref{prop:chamber} we can suppose $H\in\mathcal{C}$. The proof is then essentially the one of Lemma II.6 of \cite{OG1}: there one requires $\xi$ to be primitive, but Yoshioka noticed that the same argument goes through in the more general case of $r$ and $\xi$ prime to each other (see \cite{Y3}). This last condition is always verified: write $w=(r,\xi,a)$, which is primitive and $w^{2}=2$, i. e. $\xi^{2}=2ra+2$. Suppose that $s\in\mathbb{N}$ is such that $r=sr'$ and $\xi=s\xi'$. Then $s^{2}(\xi')^{2}=2sar'+2$. As $S$ is abelian or K3, we have $(\xi')^{2}=2l$ for some $l\in\mathbb{Z}$, so that $s(sl-ar')=1$. As $s\in\mathbb{N}$ this implies $s=1$, and we are done .\endproof

An important class of OLS-triples is given by those on elliptic K3 or abelian surfaces, as in this case we hae a priviledged class of polarizations. In order to prove that the deformation class of $\widetilde{M}_{v}(S,H)$ depends only on the rank of $v$, the strategy will be to deform the OLS-triple $(S,v,H)$ to an OLS-triple on an elliptic K3 or abelian surface with a polarization in this priviledged class. Let then $Y$ be an elliptic K3 or abelian surface such that $NS(Y)=\mathbb{Z}\cdot f\oplus\mathbb{Z}\cdot\sigma$, where $f$ is the class of a fiber and $\sigma$ is the class of a section. Let $v$ be a Mukai vector on $Y$, and recall the following definition (see \cite{OG1}):

\begin{defn}
A polarization $H$ on $Y$ is called $v-$\textit{suitable} if $H$ is in the unique $v-$chamber whose closure contains $f$.
\end{defn}

We have an easy numerical criterion to guarantee that a polarization on $Y$ is $v-$suitable:

\begin{lem}
\label{lem:vsuit}Let $Y$ be an elliptic K3 or abelian surface such that $NS(Y)=\mathbb{Z}\cdot\sigma\oplus\mathbb{Z}\cdot f$, where $\sigma$ is a section and $f$ is a fibre, and let $v=(v_{0},v_{1},v_{2})$ be a Mukai vector on $Y$ such that $v_{0}>0$. Let $H$ be a polarization such that $c_{1}(H)=\sigma+lf$ for some $l\in\mathbb{Z}$. 
\begin{enumerate}
 \item If $S$ is K3, then $H$ is $v-$suitable if $l\geq|v|+1$.
 \item If $S$ is abelian, then $H$ is $v-$suitable if $l\geq|v|$.
\end{enumerate}
\end{lem}

\proof If $S$ is a K3 surface, this is Lemma I.0.3 of \cite{OG1}. For the abelian case the proof is similar: $H$ is $v-$suitable if and only if $D\cdot H$ and $D\cdot f$ have the same sign for every $D\in W_{v}$. Notice that $D=a\sigma+bf$ for some $a,b\in\mathbb{Z}$, so that $D\cdot f=a$. Suppose $D\cdot f>0$, i. e. $a>0$. Then $D\cdot H=la+b$ and $D^{2}=2ab$. As $D^{2}\geq -|v|$, we then get $b\geq -l/2a$. Hence $$D\cdot H=la+b\geq |v|a-(|v|/2a)>0,$$and we are done.\endproof

The main result of this section is the following:

\begin{prop}
\label{prop:defor}Let $(S_{1},v_{1},H_{1})$ and $(S_{2},v_{2},H_{2})$ be two OLS-triples verifying the two following conditions:
\begin{enumerate}
 \item $S_{1}$ and $S_{2}$ are two projective K3 surfaces or two abelian surfaces;
 \item if $v_{i}=2(r_{i},\xi_{i},a_{i})$, then $r_{1}=r_{2}>0$.
\end{enumerate}
Then $\widetilde{M}_{v_{1}}(S_{1},H_{1})$ is deformation equivalent to $\widetilde{M}_{v_{2}}(S_{2},H_{2})$. In particular, Theorem \ref{thm:main} is true for $(S_{1},v_{1},H_{1})$ if and only if it is true for $(S_{2},v_{2},H_{2})$.
\end{prop}

\proof The argument we present here was first used by O'Grady in \cite{OG1}, then extended by Yoshioka in \cite{Y1}. First of all, we can always assume $\rho(S_{i})>1$ and $H_{i}$ in a $v_{i}-$chamber. Indeed, consider a non-trivial smooth, projective deformation $\mathscr{X}_{i}$ of $S_{i}$ along an open $1-$dimensional disc $\Delta$, and let $0\in\Delta$ be such that $\mathscr{X}_{i,0}\simeq S_{i}$. By the Main Theorem of \cite{O}, we know that the locus of $t\in\Delta$ such that $\rho(\mathscr{X}_{i,t})>1$ is dense in the classical topology of $\Delta$. If $\mathscr{H}_{i}\in Pic(\mathscr{X}_{i})$ is a deformation of $H_{i}$ and $\mathscr{L}_{i}\in Pic(\mathscr{X}_{i})$ is a deformation of the line bundle $L_{i}\in Pic(S_{i})$ such that $c_{1}(L_{i})=\xi_{i}$, then the triple $(\mathscr{X}_{i,t},v_{i,t},\mathscr{H}_{i,t})$ is an OLS-triple for all but a finite number of $t\in\Delta$: hence there is $t\in\Delta$ such that $\rho(\mathscr{X}_{i,t})>1$ and $(\mathscr{X}_{i,t},v_{i,t},\mathscr{H}_{i,t})$ is an OLS-triple. Applying Proposition \ref{prop:chamber} we can then suppose $H_{i}$ to be in a $v_{i}-$chamber. 

By Lemma \ref{lem:polgenxi} we can suppose the triples to be $(S_{i},v_{i},H_{i})$ where $v_{i}=2(r,c_{1}(H_{i}),a_{i})$, and if $H_{i}^{2}=2d_{i}$, then we suppose $d_{i}\gg 0$. Now, let $Y$ be a K3 (resp. abelian) surface admitting an elliptic fibration and such that $$NS(Y)=\mathbb{Z}\cdot\sigma\oplus\mathbb{Z}\cdot f,$$where $f$ is the class of a fiber, and $\sigma$ is the class of a section. For $i=1,2$, there is a smooth, connected curve $T_{i}$ and a deformation $(\mathscr{X}_{i},\mathscr{L}_{i},\mathscr{H}_{i})$ over $T_{i}$ of the OLS-triple $(S_{i},v_{i},H_{i})$ such that there is $t\in T_{i}$ with the property $(\mathscr{X}_{i,t},v_{i,t},H_{i,t})=(Y,v'_{i},H_{i}')$, where  
\begin{enumerate}
 \item $c_{1}(H'_{i})=\sigma+l_{i}f$, where $l_{i}\gg 0$.
 \item $v'_{i}=2(r,c_{1}(H'_{i}),a_{i})$.
\end{enumerate}
Let $\xi'_{i}:=c_{1}(H'_{i})$. Notice that $(v'_{1})^{2}=(v'_{2})^{2}$ and they have the same rank: hence $|v'_{1}|=|v'_{2}|$, so that by Lemma \ref{lem:vsuit} a polarization is $v'_{1}-$suitable if and only if it is $v'_{2}-$suitable. Again by Lemma \ref{lem:vsuit}, we have that $H'_{i}$ is $v'_{i}-$suitable for $i=1,2$, as $l_{i}\gg 0$. Then $H'_{1}$ and $H'_{2}$ lie in the same chamber $\mathcal{C}$ (which is clearly a $v'_{i}-$chamber for $i=1,2$). By Proposition \ref{prop:chamber} we then change to a common generic polarization $H\in\mathcal{C}$, which is $v'_{i}-$generic for $i=1,2$ by Lemma \ref{lem:chambervgen}.

As $(v'_{1})^{2}=(v'_{2})^{2}$, we have $(\xi'_{1})^{2}-2ra_{1}=(\xi'_{2})^{2}-2ra_{2}$, and as $$(\xi'_{i})^{2}=(\sigma+l_{i}f)^{2}=2l_{i}-2,$$(in the abelian case we have $\xi_{i}^{2}=2l_{i}$), we then get the equation
\begin{equation}
\label{eq:k1k2}l_{1}=l_{2}+r(a_{1}-a_{2}).
\end{equation}
Notice that $v'_{1}$ and $v'_{2}$ are then equivalent: indeed, we have $$v'_{2}\cdot ch(\mathscr{O}_{Y}((a_{1}-a_{2})f))=2(r,\sigma+l_{2}f,a_{2})\cdot(1,(a_{1}-a_{2})f,0)=$$ $$=2(r,\sigma+l_{1}f,a_{1})=v'_{1},$$where the second equality follows from equation (\ref{eq:k1k2}). By Lemma \ref{lem:mustbir}, we are then done.\endproof

\begin{oss}
\label{oss:uso}{\rm We observe that in order to relate $\widetilde{M}_{v_{1}}(S_{1},H_{1})$ and $\widetilde{M}_{v_{2}}(S_{2},H_{2})$ in the previous proof, we only used deformations of the symplectic resolutions induced by deformations of OLS-triples along a smooth, connected curve, and isomorphisms between moduli spaces given by tensorization with a line bundle.}
\end{oss}

\subsection{Proof of Theorem \ref{thm:main}}

In this section we finally prove Theorem \ref{thm:main}: the crucial facts are two lemmas due to Yoshioka \cite{Y3}. If $S$ is an abelian or projective K3 surface, write $\Delta$ for the diagonal of $S\times S$, $\mathscr{I}_{\Delta}$ for the ideal sheaf of $\Delta$ and, if $S$ is abelian, let $\mathscr{P}$ be the Poincar\'e bundle on $S\times\widehat{S}$.

\begin{lem}
\label{lem:yoshi1}(\textbf{Yoshioka}). Let $(S,v,H)$ be an OLS-triple where $v=2(0,\xi,a)$ and $H$ is in a $v-$chamber. Moreover, let $\widehat{v}:=2(a,\xi,0)$, and suppose $a\gg 0$.
\begin{enumerate}
 \item If $S$ is K3, then the Fourier-Mukai transform $\mathcal{F}:D^{b}(S)\longrightarrow D^{b}(S)$ with kernel $\mathscr{I}_{\Delta}$ sends any $H-$(semi)stable sheaf with Mukai vector $v$ to an $H-$(semi)stable sheaf with Mukai vector $\widehat{v}$. In particular, it defines an isomorphism between $\widetilde{M}_{v}(S,H)$ and $\widetilde{M}_{\widehat{v}}(S,H)$.
 \item If $S$ is abelian, then the Fourier-Mukai transform $\mathcal{F}:D^{b}(S)\longrightarrow D^{b}(\widehat{S})$ with kernel $\mathscr{P}$ sends any $H-$(semi)stable sheaf with Mukai vector $v$ to an $\widehat{H}-$(semi)stable sheaf with Mukai vector $\widehat{v}$. In particular, it defines an isomorphism between $\widetilde{K}_{v}(S,H)$ and $\widetilde{K}_{\widehat{v}}(\widehat{S},\widehat{H})$.
\end{enumerate}
\end{lem}

\proof Let $w:=v/2$ and $\widehat{w}:=\widehat{v}/2$. By Proposition 3.14 of \cite{Y3}, as $a\gg 0$ the Fourier-Mukai functor of the statement sends an $H-$stable sheaf with Mukai vector $v$ (resp. $w$) to an $H-$stable sheaf of Mukai vector $\widehat{v}$ (resp. $\widehat{w}$), hence it defines an isomorphism between $M_{v}^{s}$ and $M_{\widehat{v}}^{s}$ (resp. between $M_{w}^{s}$ and $M_{\widehat{w}}^{s}$). Notice that as $H$ is $v-$generic, then it is even $w-$generic. As $w$ and $\widehat{w}$ are primitive, we then have $M_{w}=M_{w}^{s}$ and $M_{\widehat{w}}=M_{\widehat{w}}^{s}$, so that the Fourier-Mukai transform of the statement defines an isomorphism between $M_{w}$ and $M_{\widehat{w}}$. As $M_{v}=M_{v}^{s}\cup S^{2}M_{w}$ and $M_{\widehat{v}}=M_{\widehat{v}}^{s}\cup S^{2}M_{\widehat{w}}$, we are done.\endproof

The following lemma is Theorem 3.18 of \cite{Y3}:

\begin{lem}
\label{lem:yoshi2}(\textbf{Yoshioka}). Let $S$ be an abelian or projective K3 surface such that $NS(S)=\mathbb{Z}\cdot h$, where $h=c_{1}(H)$ is ample. Let $n,r\in\mathbb{N}$ such that $a:=(n^{2}-1)/r\in\mathbb{Z}$, and suppose $n\gg 0$. 
\begin{enumerate}
 \item The Fourier-Mukai transform $\mathcal{F}:D^{b}(S)\longrightarrow D^{b}(S)$ with kernel $\mathscr{I}_{\Delta}$ sends $H-$(semi)stable sheaves with Mukai vector $2(r,nh,a)$ to $H-$(semi)stable sheaves with Mukai vector $2(a,nh,r)$. In particular, it defines an isomorphism between $\widetilde{M}_{2(r,nh,a)}(S,H)$ and $\widetilde{M}_{2(a,nh,r)}(S,H)$.
 \item The Fourier-Mukai transform $\mathcal{F}:D^{b}(S)\longrightarrow D^{b}(\widehat{S})$ with kernel $\mathscr{P}$ sends $H-$(semi)stable sheaves with Mukai vector $2(r,nh,a)$ to $\widehat{H}-$(semi)stable sheaves with Mukai vector $2(a,nh,r)$. In particular, it defines an isomorphism between $\widetilde{K}_{2(r,nh,a)}(S,H)$ and $\widetilde{K}_{2(a,nh,r)}(\widehat{S},\widehat{H})$.
\end{enumerate}
\end{lem}

We now proceed with the proof of Theorem \ref{thm:main}:
\par\bigskip
\noindent\textbf{Theorem 1.6.} \textit{Let} $(S,v,H)$ \textit{be an OLS-triple.}
\begin{enumerate}
 \item \textit{If} $S$ \textit{is projective K3, then} $\widetilde{M}_{v}(S,H)$ \textit{is irreducible symplectic and deformation equivalent to} $\widetilde{M}_{10}$\textit{.}
 \item \textit{If} $S$ \textit{is abelian, then} $\widetilde{K}_{v}(S,H)$ \textit{is irreducible symplectic and deformation equivalent to} $\widetilde{K}_{6}$\textit{.}
\end{enumerate}

\proof Let $S$ be a projective K3 surface (the proof in the abelian case is analogue), and write $v=2(r,\xi,a)$. We show that $\widetilde{M}_{v}(S,H)$ is deformation equivalent to $\widetilde{M}_{2(0,h,0)}(X,H)$, where $X$ is a surface such that $NS(X)=\mathbb{Z}\cdot h$, $h=c_{1}(H)$ is ample and $H^{2}=2$. The equivalence is obtained using deformations of the simplectic resolutions induced by deformations along smooth, connected curves of the corresponding OLS-triple, and isomorphism between moduli spaces. As a particular case is $\widetilde{M}_{10}$, we will be done.

\textit{Step 1}: suppose that $NS(S)=\mathbb{Z}\cdot h$ where $h=c_{1}(H)$ is ample and $H^{2}=2$, and suppose $v=2(0,h,a)$ where $a=2k$ for some $k\in\mathbb{Z}$. Then $\widetilde{M}_{v}(S,H)\simeq\widetilde{M}_{2(0,h,0)}(S,H)$: indeed, it is immediate to see that $v=2(0,h,0)\cdot ch(\mathscr{O}_{S}(kH))$. As tensoring with a multiple of $H$ does not change $H-$(semi)stability, we get an isomorphism $$M_{2(0,h,0)}(S,H)\longrightarrow M_{2(0,h,a)}(S,H),\,\,\,\,\,\,\,\,\,\,\,\,\,\,E\mapsto E\otimes\mathscr{O}_{S}(kH),$$and we are done.

\textit{Step 2}: suppose that $(S,v,H)$ is an OLS-triple such that $r>0$. By Proposition \ref{prop:defor} we know that $\widetilde{M}_{v}(S,H)$ is deformation equivalent to $\widetilde{M}_{2(r,nh,a)}(X,H)$, where $a=(n^2-1)/r$, for some $n\in\mathbb{Z}$. Choose $n\gg 0$ such that the corresponding $a$ is even. As $n\gg 0$, point 1 of Lemma \ref{lem:yoshi2} gives an isomorphism between $M_{2(r,nh,a)}(X,H)$ and $M_{2(a,nh,r)}(X,H)$, which is deformation equivalent to $M_{2(a,h,0)}(X,H)$ by Proposition \ref{prop:defor}. As $a\gg 0$, by point 1 of Lemma \ref{lem:yoshi1} we have an isomorphism between $M_{2(a,h,0)}(X,H)$ and $M_{2(0,h,a)}(X,H)$. As $a$ is even and $NS(X)=\mathbb{Z}\cdot h$ where $h=c_{1}(H)$ is ample and $H^{2}=2$, we have $M_{2(0,h,a)}(X,H)\simeq M_{2(0,h,0)}(X,H)$ by Step 1, and we are done.

\textit{Step 3}: suppose that $(S,v,H)$ is any OLS-triple such that $r=0$. Let $d\in\mathbb{N}$ and $v':=v\cdot ch(\mathscr{O}_{S}(dH))$. As tensoring with a multiple of $H$ does not change $H-$(semi)stability, we have an isomorphism $$M_{v}(S,H)\longrightarrow M_{v'}(S,H),\,\,\,\,\,\,\,\,\,\,\,\,\,E\mapsto E\otimes\mathscr{O}_{S}(dH).$$Notice that $$v'=2(0,\xi,a)\cdot (1,dH,d^{2}H^{2}/2)=2(0,\xi,a+dH\cdot\xi).$$As $H$ is ample and $\xi$ is effective, we have $\xi\cdot H>0$, so that choosing $d\gg 0$ we get $a+d\xi\cdot H\gg 0$. We can then suppose $v=2(0,\xi,a)$ where $a\gg 0$, hence we can suppose even that $H$ is in a $v-$chamber. By point 1 of Lemma \ref{lem:yoshi1} we have then an isomorphism between $M_{v}(S,H)$ and $M_{\widehat{v}}(S,H)$, where $\widehat{v}:=2(a,\xi,0)$. We are now in the situation of Step 2, hence we are done.\endproof 

\section{The second integral cohomology of the moduli spaces}

Let $(S,v,H)$ be an OLS-triple. In this section we define a morphism $$\lambda_{v}:v^{\perp}\longrightarrow H^{2}(M_{v},\mathbb{Z}).$$For primitive Mukai vectors $v$ with $v^{2}=0$, this was defined (using semi-universal families) first by Mukai \cite{M2}, who showed that it gives a Hodge isometry between $v^{\perp}/\mathbb{Z}\cdot v$ and $H^{2}(M_{v},\mathbb{Z})$ (in this case $M_{v}$ is a K3 surface). For $v$ primitive and $v^{2}\geq 2$, this morphism was constructed by Mukai \cite{M2}, O'Grady \cite{OG1} and Yoshioka \cite{Y1}, who showed that $\lambda_{v}$ gives a Hodge isometry between $v^{\perp}$ and $H^{2}(M_{v},\mathbb{Z})$ (the latter being a lattice with respect to the Beauville form, as it is an irreducible symplectic manifold).

In the present section we define $\lambda_{v}$ for any OLS-triple $(S,v,H)$: as in the case of primitive Mukai vectors, using a semi-universal family on $S\times M_{v}^{s}$ one defines a morphism $\lambda_{v}^{s}$ which, a priori, takes values only in $H^{2}(M_{v}^{s},\mathbb{Q})$. As $(S,v,H)$ is an OLS-triple, $M_{v}^{s}$ is an open subset which is strictly contained in $M_{v}$, hence the problem is to extend $\lambda_{v}^{s}$ to a morphism $\lambda_{v}$ which takes values in $H^{2}(M_{v},\mathbb{Z})$. To show that such an extension exists, we show first, using the Le Potier morphism, that if $\alpha\in(v^{\perp})^{1,1}$, then $\lambda_{v}^{s}(\alpha)$ extends to a unique class $\lambda_{v}(\alpha)\in H^{2}(M_{v},\mathbb{Z})$ (which is, moreover, the first Chern class of a line bundle on $M_{v}$). Using a deformation argument we extend this result to every $\alpha\in v^{\perp}$. The details of the construction are explained in section 3.2.

The main result of the section is to show that the morphism $\lambda_{v}$ is a Hodge isometry between $v^{\perp}$ and $H^{2}(M_{v},\mathbb{Z})$ (or $H^{2}(K_{v},\mathbb{Z})$ if $S$ is abelian). Before doing this, we need to show that on $H^{2}(M_{v},\mathbb{Z})$ one has a pure weight-two Hodge structure and a lattice structure: as shown in section 3.1, these are induced by the respective structures on $H^{2}(\widetilde{M}_{v},\mathbb{Z})$ (the lattice structure given by the Beauville form, as $\widetilde{M}_{v}$ is now known to be an irreducible symplectic manifold), as a consequence of the fact that the singularities of $M_{v}$ are rational. In section 3.3 we show that $\lambda_{v}$ is a Hodge isometry by following the steps of the proof of Theorem \ref{thm:main}: namely, to show that $\lambda_{v}$ is a Hodge isometry is equivalent to show the same property for the O'Grady examples; in these cases, this follows immediately by \cite{Pe}.

\subsection{Hodge and lattice structures}

In this section we show that $H^{2}(M_{v},\mathbb{Z})$ and $H^{2}(K_{v},\mathbb{Z})$ admit a pure weight-two Hodge structure and a lattive structure for any OLS-triple $(S,v,H)$. As a first step we show that they are free $\mathbb{Z}-$modules.

\begin{lem}
\label{lem:ratsing}Let $X$ be a normal, irreducible projective variety with rational singularities, and let $f:\widetilde{X}\longrightarrow X$ be a resolution of the singularities. The morphism $f^{*}:H^{2}(X,\mathbb{Z})\longrightarrow H^{2}(\widetilde{X},\mathbb{Z})$ is injective.
\end{lem}

\proof As $X$ is a normal, irreducible projective variety having rational singularities and $f:\widetilde{X}\longrightarrow X$ is a resolution of singularities, then $R^{i}f_{*}\mathscr{O}_{\widetilde{X}}=0$ for every $i>0$. Moreover, by the Zariski Main Theorem the natural morphism $\mathscr{O}_{X}\longrightarrow f_{*}\mathscr{O}_{\widetilde{X}}$ is an isomorphism, and $f_{*}\mathscr{O}^{*}_{\widetilde{X}}\simeq\mathscr{O}^{*}_{X}$. Applying the functor $Rf_{*}$ to the exponential sequence of $\widetilde{X}$ we then find $R^{1}f_{*}\mathbb{Z}=0$. Using the Leray spectral sequence $$E^{p,q}_{2}:=H^{p}(X,R^{q}f_{*}\mathbb{Z})\Longrightarrow H^{p+q}:=H^{p+q}(\widetilde{X},\mathbb{Z}),$$we have $E^{p,1}_{2}=0$ for every $p\in\mathbb{Z}$, so that the map $E^{2,0}_{2}\longrightarrow H^{2}$ is injective. But this is the map $f^{*}:H^{2}(X,\mathbb{Z})\longrightarrow H^{2}(\widetilde{X},\mathbb{Z})$, and we are done.\endproof

\begin{cor}
\label{cor:h2248}Let $(S,v,H)$ be an OLS-triple.
\begin{enumerate}
 \item If $S$ is K3, then $H^{2}(M_{v},\mathbb{Z})$ is free.
 \item If $S$ is abelian, then $H^{2}(K_{v},\mathbb{Z})$ is free.
\end{enumerate}
\end{cor}

\proof If $S$ is a K3 surface, then $M_{v}$ has rational singularities: indeed, it admits a symplectic resolution, hence the singularities are canonical, so it has rational singularities by Elkik \cite{E}. By Lemma \ref{lem:ratsing}, $$\pi_{v}^{*}:H^{2}(M_{v},\mathbb{Z})\longrightarrow H^{2}(\widetilde{M}_{v},\mathbb{Z})$$is injective. Finally, by Theorem \ref{thm:main} we know that $\widetilde{M}_{v}$ is an irreducible symplectic manifold, hence it is simply connected. This implies that $H^{2}(\widetilde{M}_{v},\mathbb{Z})$ is free, so $H^{2}(M_{v},\mathbb{Z})$ is free. The case of abelian surfaces is analogue.\endproof

\begin{oss}
{\rm By Lemma \ref{lem:ratsing} and the proof of Corollary \ref{cor:h2248}, the pull-back morphism $\pi^{*}_{v}:H^{2}(M_{v},\mathbb{Z})\longrightarrow H^{2}(\widetilde{M}_{v},\mathbb{Z})$ is an injection of mixed Hodge structures. By strict compatibility of the weight filtrations, the mixed Hodge structure on $H^{2}(M_{v},\mathbb{Z})$ is then pure of weight two. Explicitely, the pure weight-two Hodge structure on $H^{2}(M_{v},\mathbb{Z})$ is defined as follows:}
\end{oss}

\begin{defn}
Let $(S,v,H)$ be an OLS-triple where $S$ is a K3 surface. The pure weight-two Hodge structure on $H^{2}(M_{v},\mathbb{Z})$ is defined as follows: $$H^{2,0}(M_{v}):=\pi_{v}^{*}(H^{2}(M_{v},\mathbb{C}))\cap H^{2,0}(\widetilde{M}_{v}),$$ $$H^{1,1}(M_{v}):=\pi_{v}^{*}(H^{2}(M_{v},\mathbb{C}))\cap H^{1,1}(\widetilde{M}_{v}),$$ $$H^{0,2}(M_{v}):=\pi_{v}^{*}(H^{2}(M_{v},\mathbb{C}))\cap H^{0,2}(\widetilde{M}_{v}).$$Similarily, we define the pure weight-two Hodge structure on $H^{2}(K_{v},\mathbb{Z})$.
\end{defn}

We now deal with the lattice structure. Recall that if $(S,v,H)$ is an OLS-triple, then $\widetilde{M}_{v}(S,H)$ and $\widetilde{K}_{v}(S,H)$ are irreducible symplectic manifolds. This implies that on $H^{2}(\widetilde{M}_{v},\mathbb{Z})$ and $H^{2}(\widetilde{K}_{v},\mathbb{Z})$ we have a lattice structure with respect to the Beauville form. As $\pi_{v}^{*}$ is injective, this induces a lattice structure on $H^{2}(M_{v},\mathbb{Z})$ and $H^{2}(K_{v},\mathbb{Z})$ which is compatible with the Hodge structure we just defined. More explicitely, we have the following:

\begin{defn}
Let $(S,v,H)$ be an OLS-triple where $S$ is a K3 surface. We define a lattice structure on $H^{2}(M_{v},\mathbb{Z})$ in the following way: define a quadratic form $q_{v}$ on $H^{2}(M_{v},\mathbb{Z})$ by letting, for every $\alpha,\beta\in H^{2}(M_{v},\mathbb{Z})$, $$q_{v}(\alpha,\beta):=\widetilde{q}_{v}(\pi_{v}^{*}\alpha,\pi_{v}^{*}\beta),$$where $\widetilde{q}_{v}$ is the Beauville form of $\widetilde{M}_{v}(S,H)$. Similarily, we define a lattice structure on $H^{2}(K_{v},\mathbb{Z})$ for every OLS-triple $(S,v,H)$ where $S$ is an abelian surface.
\end{defn}

\subsection{Mukai-Donaldson-Le Potier morphism}

In this section we define a morphism $$\lambda_{v}:v^{\perp}\longrightarrow H^{2}(M_{v},\mathbb{Z})$$for every OLS-triple $(S,v,H)$. The strategy is the following: consider a semi-universal family $\mathscr{F}$ on $S\times M_{v}^{s}$ of similitude $\rho$. Then define $$\lambda^{s}_{v,\mathscr{F}}:H^{2*}(S,\mathbb{Z})\longrightarrow H^{2}(M^{s}_{v},\mathbb{Q}),\,\,\,\,\,\lambda^{s}_{v,\mathscr{F}}(\alpha):=\frac{1}{\rho}[p_{M*}(p_{S}^{*}(\alpha^{\vee}\cdot\sqrt{td{S}})\cdot ch(\mathscr{F}))]_{1}.$$Here, if $\alpha=(\alpha_{0},\alpha_{1},\alpha_{2})$, we define $\alpha^{\vee}:=(\alpha_{0},-\alpha_{1},\alpha_{2})$, and $p_{M}$ and $p_{S}$ are the two projections of $S\times M_{v}^{s}$ to $M_{v}^{s}$ and $S$ respectively. If $S$ is abelian, composing with the inclusion morphism $j^{s}_{v}:K^{s}_{v}\longrightarrow M^{s}_{v}$ we then get a morphism $\nu^{s}_{v,\mathscr{F}}:=j_{v}^{*}\circ\lambda^{s}_{v,\mathscr{F}}$.

Now, if $\alpha\in v^{\perp}$ and $\mathscr{F},\mathscr{F}'$ are two semi-universal families, then $\lambda^{s}_{v,\mathscr{F}}(\alpha)=\lambda^{s}_{v,\mathscr{F}'}(\alpha)$ (resp. $\nu^{s}_{v,\mathscr{F}}(\alpha)=\nu^{s}_{v,\mathscr{F}'}(\alpha)$). We have then a map $$\lambda^{s}_{v}:v^{\perp}\longrightarrow H^{2}(M^{s}_{v},\mathbb{Q}),\,\,\,\,\,\,\,\,\,({\rm resp.}\,\,\,\nu^{s}_{v}:v^{\perp}\longrightarrow H^{2}(K^{s}_{v},\mathbb{Z}))$$which does not depend on the chosen semi-universal family. The problem is to extend $\lambda_{v}^{s}$ to a morphism $$\lambda_{v}:v^{\perp}\longrightarrow H^{2}(M_{v},\mathbb{Z}),$$i. e. such that if $i_{v}:M_{v}^{s}\longrightarrow M_{v}$ is the inclusion, we have $\lambda^{s}_{v}=i_{v}^{*}\circ\lambda_{v}$.  If $S$ is abelian, and $j_{v}:K_{v}\longrightarrow M_{v}$ is the inclusion, we then get a morphism $$\nu_{v}:=j_{v}^{*}\circ\lambda_{v}:v^{\perp}\longrightarrow H^{2}(K_{v},\mathbb{Z}).$$In order to do this, we need to study the relation between $H^{2}(M_{v})$ and $H^{2}(M_{v}^{s})$. We have the following:

\begin{lem}
\label{lem:struttura}Let $(S,v,H)$ be an OLS-triple, and let $i_{v}:M_{v}^{s}\longrightarrow M_{v}$ (resp. $i_{v}:K_{v}^{s}\longrightarrow K_{v}$) be the inclusion. Then $$i_{v}^{*}:H^{2}(M_{v},\mathbb{Z})\longrightarrow H^{2}(M_{v}^{s},\mathbb{Z})$$(resp. $i_{v}^{*}:H^{2}(K^{s}_{v},\mathbb{Z})\longrightarrow H^{2}(K_{v},\mathbb{Z})$) is injective.
\end{lem}

\proof We have a commutative diagram, every row of which is exact:
\begin{equation}
\label{eq:diah}
\begin{array}{ccccc}
H^{2}(M_{v},M_{v}^{s}) & \stackrel{c}\rightarrow & H^{2}(M_{v},\mathbb{Z}) & \stackrel{i_{v}^{*}}\rightarrow & H^{2}(M_{v}^{s},\mathbb{Z})\\
\scriptstyle{f}\downarrow & & \scriptstyle{\pi_{v}^{*}}\downarrow & & \scriptstyle{\pi_{v}^{*}}\downarrow\\
H^{2}(\widetilde{M}_{v},\pi_{v}^{-1}(M_{v}^{s})) & \stackrel{\widetilde{c}}\rightarrow & H^{2}(\widetilde{M}_{v},\mathbb{Z}) & \stackrel{\widetilde{i}_{v}^{*}}\rightarrow & H^{2}(\pi_{v}^{-1}(M_{v}^{s}),\mathbb{Z})\end{array}
\end{equation}
where $\widetilde{i}_{v}$ is the restriction morphism from $\widetilde{M}_{v}$ to $\pi_{v}^{-1}(M_{v}^{s})$. As $\widetilde{M}_{v}\setminus M_{v}^{s}\simeq\widetilde{\Sigma}_{v}$, which is irreducible, we have $H^{2}(\widetilde{M}_{v},\pi_{v}^{-1}(M_{v}^{s}))\simeq\mathbb{Z}$, and $\widetilde{c}(1)=c_{1}(\widetilde{\Sigma}_{v})$. Let $\alpha\in H^{2}(M_{v},\mathbb{Z})$ be such that $i_{v}^{*}(\alpha)=0$, so that $\widetilde{i}_{v}^{*}\circ\pi_{v}^{*}(\alpha)=0$. As the second row of the diagram (\ref{eq:diah}) is exact, there is $n\in\mathbb{Z}$ such that $\pi_{v}^{*}(\alpha)=\widetilde{c}(n)=nc_{1}(\widetilde{\Sigma}_{v})$. Now, we introduce the following notation: let $\Sigma_{v}^{0}\subseteq\Sigma_{v}$ be the smooth locus of $\Sigma_{v}$. Following \cite{OG2} we know that $\pi_{v}:\pi_{v}^{-1}(\Sigma_{v}^{0})\longrightarrow\Sigma_{v}^{0}$ is a $\mathbb{P}^{1}-$bundle, whose generic fiber is then a rational curve $\delta$. As $\delta$ is contracted by $\pi_{v}$, we have $\pi_{v}^{*}(\alpha)\cdot\delta=0$. On the other hand, by adjunction the normal bundle to $\widetilde{\Sigma}_{v}$ is the canonical bundle of $\widetilde{\Sigma}_{v}$, hence it has degree $-2$ on $\delta$. In conclusion, we have $$0=\pi_{v}^{*}(\alpha)\cdot\delta = nc_{1}(\widetilde{\Sigma}_{v})\cdot\delta=-2n,$$so that $n=0$. Hence $\pi_{v}^{*}(\alpha)=0$, but as $\pi_{v}^{*}$ is injective by Lemma \ref{lem:ratsing}, we then have $\alpha=0$, so that $i_{v}^{*}$ is injective, and we are done for the K3 surface case. The proof of the abelian case is similar.\endproof

By Lemma \ref{lem:struttura}, if $\alpha\in v^{\perp}$ and $\mu_{1}(\alpha),\mu_{2}(\alpha)\in H^{2}(M_{v},\mathbb{Z})$ are such that $i_{v}^{*}(\mu_{1}(\alpha))=i_{v}^{*}(\mu_{2}(\alpha))$, then by Lemma \ref{lem:struttura} we have that $\mu_{1}(\alpha)=\mu_{2}(\alpha)$. Hence, if there is an extension of $\lambda_{v}^{s}(\alpha)$ to an element of $H^{2}(M_{v},\mathbb{Z})$, then this extension is unique, and we call it $\lambda_{v}(\alpha)$. In conclusion, the problem is only to find an extension.

In order to do so, we recall a construction due to Le Potier. Let $K_{hol}(S)$ be the holomorphic $K-$theory of $S$, and let $$vect^{\vee}:K_{hol}(S)\longrightarrow H^{2*}(S,\mathbb{Z}),\,\,\,\,\,\,\,\,\,\,\,vect^{\vee}([E]):=(v(E))^{\vee},$$where $[E]$ is the class in $K_{hol}(S)$ of a sheaf $E$ on $S$. Notice that $vect^{\vee}$ gives an isomorphism between $K_{hol}(S)$ and $\widetilde{H}^{1,1}(S)\cap H^{2*}(S,\mathbb{Z})$. Let $$(.,.):K_{hol}(S)\times K_{hol}(S)\longrightarrow\mathbb{Z},\,\,\,\,\,\,\,\,\,\,\,\,([E],[F]):=\chi(E\otimes F),$$and it is easy to see that $([E],[F])=-(v(E),v(F))$ for every sheaves $E,F$ on $S$. Let $\mathscr{E}$ be any sheaf parameterized by $M_{v}(S,H)$, $e_{v}:=[\mathscr{E}]$ and $e_{v}^{\perp}\subseteq K_{hol}(S)$ the orthogonal of $e_{v}$ with respect to $(.,.)$. Finally, let $Q_{v}$ be a Quot-scheme such that $M_{v}=Q_{v}/GL(N)$ for some $N\in\mathbb{Z}$, and let $R_{v}\subseteq Q_{v}$ be the open subset of $Q_{v}$ parameterizing $H-$semistable quotients. Let $q_{R}$ and $q_{S}$ be the projections of $S\times R_{v}$ onto $R_{v}$ ad $S$ respectively, and let $\mathscr{F}$ be a universal family on $S\times R_{v}$. Then define $$L^{R}_{v,\mathscr{F}}:K_{hol}(S)\longrightarrow Pic(R_{v}),\,\,\,\,\,L^{R}_{v,\mathscr{F}}([E]):=det(p_{R!}(p_{S}^{*}[E]\cdot[\mathscr{F}])).$$If $\mathscr{F}$ and $\mathscr{F}'$ are two universal families on $S\times R_{v}$, then $L^{R}_{v,\mathscr{F}}([E])=L^{R}_{v,\mathscr{F}'}([E])$ for every $[E]\in e_{v}^{\perp}$ (see Lemma 1.2 of \cite{LP}), so we get a morphism $$L^{R}_{v}:e_{v}^{\perp}\longrightarrow Pic(R_{v}).$$

\begin{lem}
\label{lem:lepotier}Let $(S,v,H)$ be any OLS-triple. Then for every $[E]\in e_{v}^{\perp}$ the line bundle $L^{R}_{v}([E])$ descends to a line bundle $L_{v}([E])\in Pic(M_{v})$.
\end{lem}

\proof The line bundle $L^{R}_{v}([E])$ has a natural $GL(N)-$linearization inherited from the one we have on $\mathscr{F}$. Let $[P]\in R_{v}$ be a point with closed $GL(N)-$orbit corresponding to a sheaf $F\in M_{v}$. Let $\pi:R_{v}\longrightarrow M_{v}$ be the quotient morphism, so that $\pi([P])=F$. We need to show that the action of the stabilizer $Stab([P])$ is trivial on the fiber $L^{R}_{v}([E])_{[P]}$. We know that this is the case if $F$ is $H-$stable by \cite{LP}, hence we suppose $F=(F_{1}\otimes V_{1})\oplus(F_{2}\otimes V_{2})$, where $F_{1},F_{2}$ are $H-$stable and $V_{1},V_{2}$ are vector spaces. We know that $$Stab([P])\simeq Aut(F)\simeq GL(V_{1})\times GL(V_{2}).$$Moreover, we have $$L^{R}_{v}([E])_{[P]}\simeq\bigotimes_{i=1}^{2}\bigg(det(H^{\bullet}(F_{i}\otimes E))^{dim(V_{i})}\otimes(det(V_{i}))^{\chi(F_{i}\otimes E)}\bigg),$$and the action of an element $(M_{1},M_{2})\in Stab([P])$ is simply the multiplication by $det(M_{1})^{\chi(F_{1}\otimes E)}det(M_{2})^{\chi(F_{2}\otimes E)}$. As the polarization $H$ is $v-$generic, then $v(F_{1})=v(F_{2})=v/2$: hence, as $[E]\in e_{v}^{\perp}$, we get $\chi(F_{1}\otimes E)=\chi(F_{2}\otimes E)=0$. In conclusion, the action of any element of the stabilizer is trivial, so that $L^{R}_{v}([E])$ descends to a line bundle $L_{v}([E])\in Pic(M_{v})$. \endproof

\begin{oss}
{\rm We observe that the argument of the previous lemma works for every Mukai vector $v$ and every $v-$generic polarization.}
\end{oss}

We have, in conclusion, a morphism $L_{v}:e_{v}^{\perp}\longrightarrow Pic(M_{v})$. The main result of the section is the following:

\begin{prop}
\label{prop:extension}
Let $(S,v,H)$ be any OLS-triple. Then there is a morphism $$\lambda_{v}:v^{\perp}\longrightarrow H^{2}(M_{v},\mathbb{Z})$$such that $i_{v}^{*}\circ\lambda_{v}=\lambda^{s}_{v}$.
\end{prop}

\proof By Lemma \ref{lem:lepotier} we have a morphism $L_{v}:e_{v}^{\perp}\longrightarrow Pic(M_{v})$. An application of the Grothendieck-Riemann-Roch Theorem shows that if $\alpha\in (v^{\perp})^{1,1}$ and $[E]\in e_{v}^{\perp}$ is the unique element such that $vect^{\vee}([E])=\alpha$, then $$\lambda^{s}_{v}(\alpha)=i_{v}^{*}(c_{1}(L_{v}([E])))$$(for a similar computation, see \cite{Pe}). Hence, we define $\lambda_{v}(\alpha):=c_{1}(L_{v}([E]))$, so that we finally get a morphism $$\lambda_{v}:(v^{\perp})^{1,1}\longrightarrow H^{2}(M_{v},\mathbb{Z}).$$It remains to show that we can define $\lambda_{v}$ on the whole $v^{\perp}$. To do so, we use a deformation argument. Let $$\Omega:=\{\omega\in\mathbb{P}(H^{2}(S,\mathbb{C}))\,|\,\omega\cup\omega=0,\,\,\omega\cup\overline{\omega}>0\}$$be the open subset of a quadric in $\mathbb{P}(H^{2}(S,\mathbb{C}))$ which parameterizes periods of K3 surfaces. Let $\mathcal{K}$ be the Kuranishi family of for $S$, and let $p:\mathcal{K}\longrightarrow\Omega$ be the period map. Let $t_{0}:=p(S)\in\Omega$, and let $C$ be a smooth, connected curve in $\Omega$ passing through $t_{0}$. Let $(\mathscr{X},\mathscr{L},\mathscr{H})$ be a deformation of $(S,v,H)$ along $C$: for every $t\in C$, we have that $t=p(\mathscr{X}_{t})$ is the period of $\mathscr{X}_{t}$, and notice that $v_{t}\in\widetilde{H}^{1,1}(\mathscr{X}_{t})$. Let $f:\mathscr{X}\longrightarrow C$ and let $\phi:\mathscr{M}\longrightarrow C$ be the relative moduli space of semistable sheaves with Mukai vector $v$. By Corollary \ref{cor:banaletop} there is an analytic open neighborhood $U$ of $t_{0}\in C$ such that $f^{-1}(U)$ is homeomorphic to the product $S\times U$ over $U$, and $\phi^{-1}(U)$ is homeomorphic to the product $M_{v}\times U$ over $U$. Up to shrinking $U$, we can even suppose that the local systems $Rf_{*}\mathbb{Z}$, $R^{2}f{*}\mathbb{Z}$, $R^{2}\phi_{*}\mathbb{Z}$ and $R^{2}\phi^{s}_{*}\mathbb{Z}$ are constant (hence even those with coefficients in $\mathbb{Q}$, $\mathbb{R}$ and $\mathbb{C}$ are constant). Notice that this means that we can identify $H^{2*}(\mathscr{X}_{t},\mathbb{Z})$ with $H^{2*}(S,\mathbb{Z})$ (as lattices), $H^{2}(\mathscr{X}_{t},\mathbb{C})$ with $H^{2}(S,\mathbb{C})$ (as lattices), $H^{2}(M_{v_{t}},\mathbb{Z})$ with $H^{2}(M_{v},\mathbb{Z})$ and $H^{2}(M_{v_{t}}^{s},\mathbb{Q})$ with $H^{2}(M^{s}_{v},\mathbb{Q})$. As $v$ is constant over $U$, we can then consider $\mathcal{V}\subseteq Rf_{*}\mathbb{Z}$ to be a local system such that for every $t\in T$ we have $\mathcal{V}_{t}=v_{t}^{\perp}$. As we have relative semi-universal families, then we define a morphism $$\lambda:\mathcal{V}\longrightarrow R^{2}\phi^{s}_{*}\mathbb{Q},$$such that for every $t\in U$ we have $\lambda_{t}=\lambda_{v_{t}}$. Notice that we just need to show that there is a system of generators $\alpha_{1},...,\alpha_{23}\in v^{\perp}$ such that for every $i=1,...,23$ there is $t_{i}\in U$ such that $\alpha_{i}\in(v^{\perp}_{t_{i}})^{1,1}$. Indeed, by Lemma \ref{lem:lepotier} we have that $\lambda^{s}_{v_{t_{i}}}(\alpha_{i})$ extends to a unique element $\lambda_{v_{t_{i}}}(\alpha_{i})\in H^{2}(M_{v_{t_{i}}},\mathbb{Z})$. Hence even $\lambda^{s}_{v}(\alpha_{i})$ extends to a unique $\lambda_{v}(\alpha_{i})\in H^{2}(M_{v},\mathbb{Z})$. Now, let $\alpha\in v^{\perp}$: then there are $\mu_{1},...,\mu_{23}$ such that $$\alpha=\sum_{i=1}^{23}\mu_{i}\alpha_{i}.$$But this implies that $\lambda^{s}_{v}(\alpha)$ extends to the element $$\lambda_{v}(\alpha):=\sum_{i=1}^{23}\mu_{i}\lambda_{v}(\alpha_{i})\in H^{2}(M_{v},\mathbb{Z}),$$and we are done.

We then prove that there is a system of generators $\alpha_{1},...,\alpha_{23}\in v^{\perp}$ such that for every $i=1,...,23$ there is $t_{i}\in U$ such that $\alpha_{i}\in(v^{\perp}_{t_{i}})^{1,1}$. To do so, define $V:=v^{\perp}\otimes\mathbb{C}$, which is a $23-$dimensional complex vector space, and let $\mathbb{P}:=\mathbb{P}(V)$. Notice that for every $t\in U$ we have that $t$ is the period of $\mathscr{X}_{t}$, hence it is a $(2,0)-$class on $\mathscr{X}_{t}$: as $v_{t}\in\widetilde{H}^{1,1}(\mathscr{X}_{t})$, then $t\in v_{t}^{\perp}$, so that in conclusion we have $U\subseteq\mathbb{P}$. Consider the incidence variety $I\subseteq U\times\mathbb{P}$, i. e. $$I:=\{(t,[w])\in U\times\mathbb{P}\,|\,(t,w)=0\},$$together with the two projections $g$ and $h$ to $U$ and $\mathbb{P}$ respectively. Notice that if $t\in U$, then $t$ is the period of $\mathscr{X}_{t}$, so that $$g^{-1}(t)\simeq\mathbb{P}(v_{t}^{\perp}\otimes\mathbb{C}\cap(\widetilde{H}^{2,0}(\mathscr{X}_{t})\oplus \widetilde{H}^{1,1}(\mathscr{X}_{t})))\subseteq\mathbb{P}(v_{t}^{\perp}\otimes\mathbb{C})\simeq\mathbb{P}.$$Hence this implies that $I$ is a smooth variety, and $dim(I)=dim(\mathbb{P})$.

Now, let us suppose that there is a point $w\in V\cap\widetilde{H}^{1,1}(S)\cap H^{2*}(S,\mathbb{R})$ such that $[w]\in im(h)$ and it admits an open analytic neighborhood $\mathcal{U}$ which is contained in $im(h)$. Then then implies that there are $\alpha_{1},...,\alpha_{23}\in v^{\perp}$ which form a system of generators, and such that $[\alpha_{i}]\in\mathcal{U}$ for every $i=1,...,23$. Then $[\alpha_{i}]\in im(h)$, which implies that for every $i=1,...,23$ there is $t_{i}\in U$ such that $(t_{i},[\alpha_{i}])\in I$, i. e. $\alpha_{i}$ is orthogonal to $t_{i}$ with respect to the Mukai pairing on $\widetilde{H}(\mathscr{X}_{t_{i}})$. As $t_{i}$ is the period of $\mathscr{X}_{t_{i}}$, this means that $\alpha_{i}\in\widetilde{H}^{2,0}(\mathscr{X}_{t_{i}})\oplus\widetilde{H}^{1,1}(\mathscr{X}_{t_{i}})$. Finally, recall that $\alpha_{i}\in v^{\perp}$, so that $\alpha_{i}\in v_{t_{i}}^{\perp}$, hence $\alpha_{i}=\overline{\alpha_{i}}$. But this implies that $\alpha_{i}\in\widetilde{H}^{1,1}(\mathscr{X}_{t_{i}})$. In conclusion, we have $\alpha_{i}\in(v_{t_{i}}^{\perp})^{1,1}$. 

In conclusion, we just need to prove that there is $w\in V\cap\widetilde{H}^{1,1}(S)\cap H^{2*}(S,\mathbb{R})$ such that $[w]\in im(h)$ and which admits an open analytic neighborhood $\mathcal{U}$ contained in $im(h)$. This is an immediate consequence of the following:

\begin{lem}
\label{lem:esiste}
There is $w\in V\cap\widetilde{H}^{1,1}(S)\cap H^{2*}(S,\mathbb{R})$ such that the map $$dh_{(t_{0},[w])}:T_{(t_{0},[w])}I\longrightarrow T_{[w]}\mathbb{P}$$is an isomorphism.
\end{lem}

\proof Let $T$ be the tangent line at $U$ in $t_{0}$. As $U\subseteq\mathbb{P}$, then there is a $2-$dimensional linear subspace $W\subseteq V$ such that $\mathbb{P}(W)=T$. Notice that as $t_{0}$ is the period of $S$, then  $T\subseteq\mathbb{P}(V\cap(\widetilde{H}^{2,0}(S)\oplus\widetilde{H}^{1,1}(S)))$, so that in conclusion $W\subseteq V\cap(\widetilde{H}^{2,0}(S)\oplus\widetilde{H}^{1,1}(S))$. Let $W^{\perp}$ be the orthogonal to $W$ with respect to the Mukai pairing on $V$. As $dim(W)=2$, then $dim(W^{\perp})=21$. If for every $w\in V\cap\widetilde{H}^{1,1}(S)$ we have that $w\in W^{\perp}$, then $W^{\perp}=V\cap\widetilde{H}^{1,1}(S)$ (as $dim(V\cap\widetilde{H}^{1,1}(S))=21$). But $(V\cap\widetilde{H}^{1,1}(S))^{\perp}=V\cap(\widetilde{H}^{2,0}(S)\oplus\widetilde{H}^{0,2}(S))$, hence $W=V\cap(\widetilde{H}^{2,0}(S)\oplus\widetilde{H}^{0,2}(S))$, which is not possible. In conclusion, there is a $w\in V\cap\widetilde{H}^{1,1}(S)$ such that $w\notin W^{\perp}$. Moreover, we can suppose that $w\in H^{2*}(S,\mathbb{R})$: indeed, as $H^{1,1}(S)$ is defined over $\mathbb{R}$, then if $V\cap\widetilde{H}^{1,1}(S)\cap H^{2*}(S,\mathbb{R})\subseteq W^{\perp}$, then we have $V\cap\widetilde{H}^{1,1}(S)\subseteq W^{\perp}$, which is not possible as we have just seen.

Consider $[w]\in\mathbb{P}$, and let $M\subseteq V$ be the line such that $\mathbb{P}(M)=[w]$, so that $M$ is not contained in $W^{\perp}$. Moreover, let $L\subseteq V$ be the line such that $\mathbb{P}(L)=t_{0}$. Recall that $$T_{[w]}\mathbb{P}\simeq Hom(M,V/M),$$ $$T_{(t_{0},[w])}(U\times\mathbb{P})\simeq Hom(L,W/L)\times Hom(M,V/M).$$By definition of $I$, we then have $$T_{(t_{0},[w])}I=\{(\phi,\psi)\in T_{(t_{0},[w])}(U\times\mathbb{P})\,|\,\,(\phi(l),m)-(l,\psi(m))=0\},$$where the equation is true for every $l\in L$ and every $m\in M$. Moreover, we have $$dh_{(t_{0},[w])}:T_{(t_{0},[w])}I\longrightarrow T_{[w]}\mathbb{P},\,\,\,\,\,\,\,\,\,\,\,dh_{(t_{0},[w])}(\phi,\psi)=\psi.$$As $I$ is smooth and $dim(I)=dim(\mathbb{P})$, in order to show that $dh_{(t_{0},[w])}$ is an isomorphism we just need to show that it is surjective. Consider then $\psi\in Hom(M,V/M)$: as $M$ is not contained in $W^{\perp}$, for every $l\in L$ there is an element $\phi(l)\in W$ such that $$(\psi(m),l)=-(m,\phi(l))$$for every $m\in M$. But this defines an element $\phi\in Hom(L,W/L)$ such that $(\phi,\psi)\in T_{(t_{0},[w])}I$ and $dh_{(t_{0},[w])}(\phi,\psi)=\psi$, and we are done.\endproof
\endproof

\subsection{Proof of Theorem \ref{thm:mainl}}

The aim of this section is to prove the following:
\par\bigskip
\noindent\textbf{Theorem 1.7.} \textit{Let} $(S,v,H)$ \textit{be an OLS-triple.}
\begin{enumerate}
 \item \textit{If} $S$ \textit{is K3, then} $\lambda_{v}:v^{\perp}\longrightarrow H^{2}(M_{v},\mathbb{Z})$ \textit{is a Hodge isometry.}
 \item \textit{If} $S$ \textit{is abelian, then} $\nu_{v}:v^{\perp}\longrightarrow H^{2}(K_{v},\mathbb{Z})$ \textit{is a Hodge isometry.}
\end{enumerate}

\proof Let $(S,v,H)$ be an OLS-triple. We need to show the three following properties:
\begin{enumerate}
 \item $\lambda_{v}$ (resp. $\nu_{v}$) is an isomorphism of $\mathbb{Z}-$modules;
 \item $\lambda_{v}$ (resp. $\nu_{v}$) is an isometry;
 \item $\lambda_{v}$ (resp. $\nu_{v}$) is a Hodge morphism.
\end{enumerate}
We introduce the following notations:
 $$\widetilde{\lambda}_{v}:=\pi_{v}^{*}\circ\lambda_{v}:v^{\perp}\longrightarrow H^{2}(\widetilde{M}_{v},\mathbb{Z}),\,\,\,\,\,\,\,\,\,\,\widetilde{\nu}_{v}:=\pi_{v}^{*}\circ\nu_{v}:v^{\perp}\longrightarrow H^{2}(\widetilde{K}_{v},\mathbb{Z}).$$
 
\textit{Step 1}. If $S$ is an abelian or projective K3 surface such that $NS(S)=\mathbb{Z}\cdot h$, where $h=c_{1}(H)$ is ample and $h^{2}=2$, and $v=(2,0,-2)$, then $\lambda_{v}$ and $\nu_{v}$ are Hodge isometries: this is proved in \cite{Pe}.

\textit{Step 2}. Let $(S,v,H)$ be an OLS-triple. We show that $\lambda_{v}$ is an isomorphism of $\mathbb{Z}-$modules. Following the proof of Theorem \ref{thm:main}, we reduce to the case of Step 1: the only transformations we use are deformations of the moduli spaces induced by deformations of the corresponding OLS-triple along a smooth, connected curve, and isomorphisms between moduli spaces induced by some Fourier-Mukai transforms. Deforming an OLS-triple along a smooth, connected curve does not change the $\mathbb{Z}-$module structures of $v^{\perp}$ and of $H^{2}(M_{v},\mathbb{Z})$; for the isomorphism induced by the Fourier-Mukai tranform we have the following:

\begin{lem}
\label{lem:fourmuk}Let $(S,v,H)$ be an OLS-triple. 
\begin{enumerate}
 \item If $S$ is K3, let $\mathscr{P}\in D^{b}(S\times S)$ and $\mathcal{F}_{\mathscr{P}}:D^{b}(S)\longrightarrow D^{b}(S)$ the Fourier-Mukai transform with kernel $\mathscr{P}$. Moreover, let $\phi_{\mathscr{P}}$ be the morphism induced by $\mathcal{F}_{\mathscr{P}}$ in cohomology, and let $v':=\phi_{\mathscr{P}}(v)$. If $\mathcal{F}_{\mathscr{P}}$ is an equivalence and it induces an isomorphism $f_{\mathscr{P}}:M_{v'}(S,H)\longrightarrow M_{v}(S,H)$, 
then $\lambda_{v}$ is an isomorphism if and only if $\lambda_{v'}$ is an isomorphism.
 \item If $S$ is abelian, let $\mathscr{P}\in D^{b}(S\times\widehat{S})$ and $\mathcal{F}_{\mathscr{P}}:D^{b}(S)\longrightarrow D^{b}(\widehat{S})$ the Fourier-Mukai transform with kernel $\mathscr{P}$. Moreover, let $\phi_{\mathscr{P}}$ be the morphism induced by $\mathcal{F}_{\mathscr{P}}$ in cohomology, and let $v':=\phi_{\mathscr{P}}(v)$. If $\mathcal{F}_{\mathscr{P}}$ is an equivalence and it induces an isomorphism $f_{\mathscr{P}}:K_{v'}(\widehat{S},\widehat{H})\longrightarrow K_{v}(S,H)$, 
then $\nu_{v}$ is an isomorphism if and only if $\nu_{v'}$ is an isomorphism.
\end{enumerate}
\end{lem}

\proof We show the first point, as the second is similar. We show that the diagram
\begin{equation}
\label{eq:diagcomm}
\begin{array}{ccc}v^{\perp} & \stackrel{\phi_{\mathscr{P}}}\longrightarrow & (v')^{\perp}\\ \scriptstyle{\lambda_{v}}\downarrow & & \downarrow\scriptstyle{\lambda_{v}}\\ H^{2}(M_{v},\mathbb{Z}) & \stackrel{f_{\mathscr{P}}^{*}}\longrightarrow & H^{2}(M_{v'},\mathbb{Z})\end{array}
\end{equation}
is commutative. By the construction of $\lambda_{v}$ and $\lambda_{v'}$, this is true if the following diagram
\begin{equation}
\label{eq:diagcomms}
\begin{array}{ccc}v^{\perp} & \stackrel{\phi_{\mathscr{P}}}\longrightarrow & (v')^{\perp}\\ \scriptstyle{\lambda^{s}_{v}}\downarrow & & \downarrow\scriptstyle{\lambda^{s}_{v}}\\ H^{2}(M^{s}_{v},\mathbb{Q}) & \stackrel{f_{\mathscr{P}}^{*}}\longrightarrow & H^{2}(M^{s}_{v'},\mathbb{Q})\end{array}
\end{equation}
is commutative, and this is shown to be true by standard computations (see for example Proposition 2.4 of \cite{Y2}). As $\phi_{\mathscr{P}}$ and $f_{\mathscr{P}}^{*}$ are isomorphisms, then $\lambda_{v}$ is an isomorphism if and only if $\lambda_{v'}$ is, and we are done.\endproof

In conclusion, we reduce to the case of Step 1, so that $\lambda_{v}$ is an isomorphism of $\mathbb{Z}-$modules for every OLS-triple $(S,v,H)$.

\textit{Step 3}. We prove now that $\lambda_{v}$ is an isometry between $v^{\perp}$ (which has a lattice structure given by the Mukai pairing) and $H^{2}(M_{v},\mathbb{Z})$ (which has lattice structure as seen in section 3.1). Again, we reduce to the case of Step 1 following the proof of Theorem \ref{thm:main}: as in the previous step, the only transformations we use are deformations of the moduli spaces induced by deformations of the corresponding OLS-triple along a smooth, connected curve, and isomorphisms between moduli spaces induced by some Fourier-Mukai transforms. As the lattice structures of $v^{\perp}$ and of $H^{2}(M_{v},\mathbb{Z})$ are locally constant along a smooth, connected curve, we just need to analyze the isomorphism induced by Fourier-Mukai transforms. We have the following:

\begin{lem}
\label{lem:isof}Let $(S_{1},v_{1},H_{1})$ and $(S_{2},v_{2},H_{2})$ be two OLS-triples. 
\begin{enumerate}
 \item If $S_{1}$ and $S_{2}$ are K3 and there is an isomorphism $f:M_{v_{1}}\longrightarrow M_{v_{2}}$, then the morphism $f^{*}:H^{2}(M_{v_{2}},\mathbb{Z})\longrightarrow H^{2}(M_{v_{1}},\mathbb{Z})$ is an isometry.
 \item If $S_{1}$ and $S_{2}$ are abelian and there is an isomorphism $f:K_{v_{1}}\longrightarrow K_{v_{2}}$, then the morphism $f^{*}:H^{2}(K_{v_{2}},\mathbb{Z})\longrightarrow H^{2}(K_{v_{1}},\mathbb{Z})$ is an isometry.
\end{enumerate}
\end{lem}

\proof We prove the first point, as the second is similar. We have a commutative diagram $$\begin{array}{ccc}H^{2}(M_{v_{2}},\mathbb{Z}) & \stackrel{f^{*}}\longrightarrow & H^{2}(M_{v_{1}},\mathbb{Z})\\ \scriptstyle{\pi_{v_{2}}^{*}}\downarrow & & \downarrow\scriptstyle{\pi_{v_{1}}^{*}}\\ H^{2}(\widetilde{M}_{v_{2}},\mathbb{Z}) & \stackrel{\widetilde{f}^{*}}\longrightarrow & H^{2}(\widetilde{M}_{v_{1}},\mathbb{Z}).\end{array}$$By hypothesis, we have that $f$ is an isomorphism. Moreover, $\pi_{v_{1}}^{*}$ and $\pi_{v_{2}}^{*}$ are morphisms of lattices, and $\widetilde{f}^{*}$ is an isometry by \cite{OG1}. Hence $f^{*}$ is an isometry, and we are done.\endproof

In conclusion, we reduce to the case of Step 1, so that $\lambda_{v}$ is an isometry for every OLS-triple $(S,v,H)$.

\textit{Step 4}. We show that $\lambda_{v}$ is a Hodge morphism. Suppose in the following that $S$ is a K3 surface (the proof for $S$ abelian is analogue). To show that $\lambda_{v}$ is a Hodge morphism is equivalent to show that $\widetilde{\lambda}_{v}$ is a Hodge morphism. Notice that as $\lambda_{v}$ is an isometry by Step 3, then $\widetilde{\lambda}_{v}$ is an isometry onto its image. Recall that $\lambda_{v}$ is defined as an extension of the morphism $$\lambda_{v}^{s}:v^{\perp}\longrightarrow H^{2}(M_{v}^{s},\mathbb{Q}),\,\,\,\,\,\,\,\,\,\lambda_{v}^{s}(\alpha)=\frac{1}{\rho}[p_{*}(q^{*}(\alpha^{\vee}\cdot\sqrt{td(S)})\cdot ch(\mathscr{F})]_{1},$$where $\rho$ is the similitude of $\mathscr{F}$ and $p,q$ are the two projections of $S\times M_{v}^{s}$ onto $M_{v}^{s}$ and $S$ respectively. As $ch(\mathscr{F})\in H^{2*}(S\times M_{v}^{s},\mathbb{Q})$, and as $M_{v}^{s}$ is (up to isomorphism) an open subset of $\widetilde{M}_{v}$, taking the closure of the cycles $ch_{i}(\mathscr{F})$ we get a class $c\in H^{2*}(S\times\widetilde{M}_{v},\mathbb{Q})$, whose component $c_{i}\in H^{2i}(S\times\widetilde{M}_{v},\mathbb{Q})$ represents a $(i,i)-$class. Let $\widetilde{p}$ and $\widetilde{q}$ be the projections of $S\times\widetilde{M}_{v}$ onto $\widetilde{M}_{v}$ and $S$ respectively, and consider the morphism $$\widetilde{\mu}_{v}:v^{\perp}\longrightarrow H^{2}(\widetilde{M}_{v},\mathbb{Q}),\,\,\,\,\,\,\,\,\,\,\,\,\widetilde{\mu}_{v}(\alpha):=\frac{1}{\rho}[\widetilde{p}_{*}(\widetilde{q}^{*}(\alpha^{\vee}\cdot\sqrt{td(S)})\cdot c)]_{1}.$$On $v^{\perp}$ and $H^{2}(\widetilde{M}_{v},\mathbb{Q})$ we have pure weight-two Hodge structures, and $\widetilde{\mu}_{v}$ is a Hodge morphism. Now, a priori the morphisms $\widetilde{\lambda}_{v}$ and $\widetilde{\mu}_{v}$ are not equal, but we have the following:

\begin{lem}
\label{lem:ugu}
For every $\omega\in(v^{\perp})^{2,0}$ we have  $\widetilde{\lambda}_{v}(\omega)=\widetilde{\mu}_{v}(\omega)$. 
\end{lem}

\proof Let $\widetilde{i}_{v}:\pi_{v}^{-1}(M_{v}^{s})\longrightarrow\widetilde{M}_{v}$ be the inclusion. By the very definition of $\widetilde{\lambda}_{v}$ and $\widetilde{\mu}_{v}$, for every $\omega\in(v^{\perp})^{2,0}$ we have $\widetilde{i}_{v}^{*}(\widetilde{\lambda}_{v}(\omega))=\widetilde{i}_{v}^{*}(\widetilde{\mu}_{v}(\omega))$. Moreover, the kernel of the morphism $\widetilde{i}_{v}^{*}:H^{2}(\widetilde{M}_{v},\mathbb{Q})\longrightarrow H^{2}(\pi_{v}^{-1}(M_{v}^{s}),\mathbb{Q})$ is $\mathbb{Q}\cdot c_{1}(\widetilde{\Sigma}_{v})$, so that $\widetilde{\lambda}_{v}(\omega)-\widetilde{\mu}_{v}(\omega)=qc_{1}(\widetilde{\Sigma}_{v})$ for some $q\in\mathbb{Q}$. But $$(\widetilde{\lambda}_{v}(\omega)-\widetilde{\mu}_{v}(\omega))\cdot\delta=0,$$so that $q=0$, and $\widetilde{\lambda}_{v}(\omega)=\widetilde{\mu}_{v}(\omega)$ for every $\omega\in(v^{\perp})^{2,0}$.\endproof 

As $\widetilde{\mu}_{v}$ is a Hodge morphism, by Lemma \ref{lem:ugu} we have that $\widetilde{\lambda}_{v}(\omega)\in H^{2,0}(\widetilde{M}_{v})$ for every $\omega\in(v^{\perp})^{2,0}$. Now, consider $\alpha\in v^{\perp}\cap (\widetilde{H}^{2,0}(S)\oplus\widetilde{H}^{1,1}(S))$. Then, for every $\omega\in(v^{\perp})^{2,0}$ we have $(\alpha,\omega)=0$, where $(.,.)$ is the Mukai pairing on $v^{\perp}$. As $\widetilde{\lambda}_{v}$ is an isometry by assumption, we have then $$\widetilde{q}_{v}(\widetilde{\lambda}_{v}(\alpha),\widetilde{\lambda}_{v}(\omega))=0,$$where $\widetilde{q}_{v}$ is the Beauville form of the irreducible symplectic manifold $\widetilde{M}_{v}$. But this implies that $\widetilde{\lambda}_{v}(\alpha)\in H^{2,0}(\widetilde{M}_{v})\oplus H^{1,1}(\widetilde{M}_{v})$ for every $\alpha\in v^{\perp}\cap(\widetilde{H}^{2,0}(S)\oplus\widetilde{H}^{1,1}(S))$. In conclusion, we see that $\widetilde{\lambda}_{v}$ respects the Hodge filtrations, hence it is a Hodge morphism, and we are done.\endproof

\subsection*{Aknowledgements}
We would like to thank Kieran O'Grady and Christoph Sorger for helpful conversation and disponibility. We thank Manfred Lehn, Yasunari Nagai, Markus Zowislok for useful conversation, and Chiara Camere for having signaled some mistakes in a preliminary draft of this work. The first author was supported by the SFB/TR 45 'Periods, Moduli spaces and Arithmetic of Algebraic Varieties' of the DFG (German Research Fundation).

\par\bigskip
\par\bigskip
Arvid Perego, Fachbereich Physik, Mathematik u. Informatik, Johannes Gutenberg Universit\"at Mainz, D-55099 Mainz, Germany.

\textit{E-mail address:} \texttt{perego@mathematik.uni-mainz.de}
\par\bigskip
Antonio Rapagnetta, Dipartimento di Matematica dell'Università di Roma II - Tor Vergata, 00133 Roma, Italy.

\textit{E-mail address:} \texttt{rapagnet@mat.uniroma2.it} 

\end{document}